\documentclass[a4paper]{article}
\usepackage{a4wide}
\usepackage[utf8]{inputenc}
\usepackage[ngerman, english]{babel}
\usepackage{amsmath,amsthm,amssymb,amsfonts}
\usepackage{mathrsfs}
\usepackage{ifthen}
\usepackage{enumerate}
\usepackage{comment}
\usepackage{subcaption}
\usepackage{hyphsubst}
\usepackage{xcolor}
\usepackage{calc}
\usepackage{graphicx}
\usepackage{hyperref}
\usepackage[toc,page]{appendix}
\usepackage{cleveref}
\usepackage{leftindex}
\usepackage{pgfplots}
\pgfplotsset{compat=1.5.1}
\usepackage[bbgreekl]{mathbbol}
\DeclareSymbolFontAlphabet{\mathbbm}{bbold}
\DeclareSymbolFontAlphabet{\mathbb}{AMSb}

\makeatletter
\def\@hspace#1{\begingroup\setlength\dimen@{#1}\hskip\dimen@\endgroup}
\makeatother

\newtheorem{theorem}{Theorem}[section]
\newtheorem{lemma}[theorem]{Lemma}

\newtheorem{corollary}[theorem]{Corollary}

\newcommand{\uten}{u}
\newcommand{\ux}{v}
\newcommand{\uy}{w}

\def\tf{^\dagger}
\def\Lspace{\calL}
\def\ses{a}
\def\op{\calA}
\def\opL{\tilde \calA}

\def\matAm{\mathbf{A}^{(m)}}
\def\matAmL{{\mathbf{\tilde A}}^{(m)}}
\def\dom{\Omega}
\def\domp{\dom_+}
\def\domm{\dom_-}
\def\dompm{\dom_\pm}
\def\domL{\tilde\Omega}
\def\domLp{\domL_+}
\def\domLm{\domL_-}

\def\hxp{h_+}
\def\hxm{h_-}
\def\hxpm{h_\pm}
\def\Nxp{N_{+}}
\def\Nxm{N_{-}}
\def\Nxpm{N_{\pm}}
\def\hy{h_{y}}
\def\Ny{M}
\def\Vx{V_{\hxpm}}
\def\VxL{\tilde V_{\hxpm}}
\def\Vy{W_{\hy}}
\newcommand{\cl}[1]{\ol{#1}^{\mathrm{cl}_{H^1_0(\dom)}}}
\def\Vsemi{\cl{\Vx\otimes H^1_0(0,\pi)}}
\def\Vfull{\Vx\otimes\Vy}

\def\VsemiL{\VxL\otimes H^1_0(0,\pi)}
\def\VfullL{\VxL\otimes\Vy}

\def\XL{\tilde X}
\def\contrast{\kappa}
\def\ratio{r}

\def\ratioy{r_y}
\def\gridx{x}
\def\gridxL{\tilde x}
\def\gridy{y}
\def\hfx{\phi}
\def\hfxL{\tilde\phi}
\def\hfy{\psi}
\def\bfy{\theta}
\def\bfydis{\hat\theta}
\def\ev{\lambda}
\def\evtmp{\tau}

\def\evdism{\hat\ev_{m,\ratioy,\hxm}}
\def\mf{\frakf}

\def\mh{\frakh}
\def\mz{\frakz}
\def\mq{\frakq}
\def\mj{\frakj}
\def\roott{t_{\contrast,\ratio}}

\def\roots{s_{\contrast,\ratio,\ratioy}}
\def\ampm{a^{(m)}_\pm}
\def\amp{a^{(m)}_+}
\def\amm{a^{(m)}_-}
\def\bmpm{b^{(m)}_\pm}
\def\bmp{b^{(m)}_+}
\def\bmm{b^{(m)}_-}
\def\mumpm{\mu_{m,\pm}}
\def\mump{\mu_{m,+}}
\def\mumm{\mu_{m,-}}
\def\mumpmI{\mu_{m,1,\pm}}
\def\mumpI{\mu_{m,1,+}}
\def\mummI{\mu_{m,1,-}}
\def\mumpmII{\mu_{m,2,\pm}}
\def\mumpII{\mu_{m,2,+}}
\def\mummII{\mu_{m,2,-}}
\def\vm{v_m}
\def\vmL{\tilde v_m}
\def\dm{d_m}
\def\dmL{\tilde d_m}
\def\betamL{\tilde \beta^{(m)}}
\def\matB{\mathbf{B}^{(\evtmp)}}
\def\adis{a}
\def\bdis{b}
\def\ampmhat{\hat a^{(m)}_\pm}
\def\amphat{\hat a^{(m)}_+}
\def\ammhat{\hat a^{(m)}_-}
\def\bmpmhat{\hat b^{(m)}_\pm}
\def\bmphat{\hat b^{(m)}_+}
\def\bmmhat{\hat b^{(m)}_-}
\def\mumpmIhat{\hat\mu_{m,1,\pm}}
\def\mumpIhat{\hat\mu_{m,1,+}}
\def\mummIhat{\hat\mu_{m,1,-}}
\def\mumpmIIhat{\hat\mu_{m,2,\pm}}
\def\mumpIIhat{\hat\mu_{m,2,+}}
\def\mummIIhat{\hat\mu_{m,2,-}}
\def\mumpmhat{\hat\mu_{m,\pm}}
\def\mumphat{\hat\mu_{m,+}}
\def\mummhat{\hat\mu_{m,-}}
\def\dmhat{\hat d_m}
\def\dmhatL{\hat{\tilde{d}}_m}
\def\vmhat{\hat v_m}
\def\vmhatL{\hat{\tilde{v}}_m}
\def\Vxzhat{\hat V_{\hxpm}^0}
\def\VxzhatL{{\hat{\tilde{V}}}_{\hxpm}^{\substack{~\\\vspace{-2.5mm}\hspace{0.2mm}0}}}

\newcommand{\ol}[1]{\overline{#1}}

\newcommand{\spl}{\langle}
\newcommand{\spr}{\rangle}
\newcommand{\bpm}{\begin{pmatrix}}
\newcommand{\epm}{\end{pmatrix}}

\renewcommand{\div}{\operatorname{div}}

\newcommand{\setN}{\mathbb{N}}

\newcommand{\setR}{\mathbb{R}}

\newcommand{\setZ}{\mathbb{Z}}

\newcommand{\calA}{\mathcal{A}}

\newcommand{\calL}{\mathcal{L}}

\newcommand{\frakf}{\mathfrak{f}}

\newcommand{\frakh}{\mathfrak{h}}

\newcommand{\frakj}{\mathfrak{j}}

\newcommand{\frakq}{\mathfrak{q}}

\newcommand{\frakz}{\mathfrak{z}}

\definecolor{brickred}{rgb}{0.8, 0.25, 0.33}
\definecolor{bostonuniversityred}{rgb}{0.8, 0.0, 0.0}
\definecolor{cornellred}{rgb}{0.7, 0.11, 0.11}
\definecolor{corn}{rgb}{0.98, 0.93, 0.36}
\definecolor{schoolbusyellow}{rgb}{1.0, 0.85, 0.0}
\definecolor{TUblue}{rgb}{0,102,153}
\colorlet{TUbluelight}{TUblue!30!white}

\usepackage{fancyhdr}
\usepackage[bbgreekl]{mathbbol}
\usepackage{authblk}
\author[1]{Martin Halla}
\author[2]{Florian Oberender}
\affil{Institut für Angewandte und Numerische Mathematik,
Karlsruher Institut für Technologie}
\affil[2]{Institut f\"ur Numerische und Angewandte Mathematik, Georg-August Universität Göttingen}
\title{On the instabilities of naive FEM discretizations for PDEs with sign-changing coefficients%
\footnote{The first author acknowledges funding from Deutsche Forschungsgemeinschaft (DFG, German Research Foundation), projects 541433971 and 258734477 – SFB 1173 and that part of this work was conducted at the Johann Radon Institute for Computational and Applied Mathematics.

The second author acknowledges support from DFG, CRC 1456 project 432680300.}}

\begin{document}

\maketitle
\begin{abstract} 
\noindent
We consider a scalar diffusion equation with a sign-changing coefficient in its principle part.
The well-posedness of such problems has already been studied extensively provided that the contrast of the coefficient is non-critical.
Furthermore, many different approaches have been proposed to construct stable discretizations thereof, because naive finite element discretizations are expected to be non-reliable in general.
However, no explicit example proving the actual instability is known and numerical experiments often do not manifest instabilities in a conclusive manner.
To this end we construct an explicit example with a broad family of meshes for which we prove that the corresponding naive finite element discretizations are unstable.
On the other hand, we also provide a broad family of (non-symmetric) meshes for which we prove that the discretizations are stable.
Together, these two findings explain the results observed in numerical experiments.
\\

\\
\noindent%
\end{abstract}


\section{Introduction}\label{sec:introduction}
In this article we consider diffusion equations $-\div(\sigma\nabla u)=f$ with a sign-changing coefficient $\sigma$, i.e., the domain $\Omega$ admits a decomposition in $\Omega_\pm$ for which $\pm\sigma|_{\Omega_\pm}>0$.
Such equations occur, e.g., for fully homogenized meta materials and their reliable simulation is essential for the development of technical devices, e.g., to control sound \cite{CummerChristensenAlu16} and for cloaking \cite{GreenleafKurylevLassasUhlmann09}.
The well-posedness of problems with sign-changing coefficients has been studied extensively by means of the T-coercivity technique \cite{BonnetBDChesnelCiarlet12,BonnetBDChesnelCiarlet14a,BonnetBDChesnelCiarlet14b} and is known to depent on the contrast of $\sigma$ and the smoothness/geometry of the interface $\ol{\Omega_+}\cap\ol{\Omega_-}$.
An alternative approach to analyze such PDEs has been investigated in \cite{Nguyen16,NguyenSil20} by means of the limiting absorption principle.
The stability of convenient finite element discretizations is known only for sufficiently large contrasts \cite{BonnetBDCiarletZwoelf10} and therefore a variety of approaches to construct stable approximations have been explored, including
locally symmetric meshes
\cite{BonnetBDCarvalhoCiarlet18,Halla23SC},
optimization based methods \cite{AbdulleHuberLemaire17,CiarletLassounonRihani23,AbdulleLemaire24,ChaabanCiarletRihani25},
boundary element methods \cite{Unger21},
weakly coercive reformulations \cite{HallaHohageOberender24}
and
primal-dual stabilizations \cite{BurmanPreussErn24}.
However, in contrast to this extensive research on the development of stable discretizations the question if those specialized methods are actually necessary has received much less attention.
Indeed, for reasonably small contrasts the error curves (for decreasing mesh sizes) of naive FEMs generally do not look reliable, but still decrease often with a saw tooth like profile \cite{BurmanPreussErn24}.
At other test cases, it can even be hard to trigger some anomalies at all \cite{HallaHohageOberender24}.
Actually, the analysis of \cite{BonnetBDCarvalhoCiarlet18}
suggest that meshes being ``almost locally symmetric'' can be expected to yield stable results, and a mesh generator might produce such meshes without further due.
However, without any quantification it is hard to obtain decisive conclusions from this observation.\\

To study such questions we construct in this article an explicit example with a piece-wise constant coefficient ($\sigma_\pm:=\sigma|_{\Omega_\pm}$) and a discretization by nodal finite elments with uniform rectangular grids in $\Omega_\pm$.
Thereby the ratio $\ratio$ of the mesh sizes in $\Omega_\pm$ will play an important role in the analysis and acts in an inverse manner to the contrast $\contrast=\sigma_+/\sigma_-$, i.e., $\ratio\contrast$ will be a crucial quantity.
We prove that depending on the parameter range that all considered discretizations are either stable or unstable.
For an unstructured mesh we expect that either case can be dominant, which explains the unconclusive observations in numerical experiments.\\

The remainder of the manuscript is structured as follows.
In \Cref{sec:setting} we specify the two considered problems and their discretization.
In particular, we consider one problem on an unbouded domain and a second problem on a bounded domain, where the first can be seen as an preperational step to the second.
In \Cref{sec:discretization} we conduct our stability analysis with \Cref{thm:full,thm:full_L} as our main results.
In \Cref{sec:num} we present computational examples to confirm our theoretical results.

\section{Notation and setting}\label{sec:setting}
We consider all vector spaces over $\setR$ and denote scalar and vectorial $L^2$-scalar products over a domain $D\subset\setR^l,l=1,2$ as $\langle\cdot,\cdot\rangle_D$.
Let $\setN:=\{1,2,\dots\}$, $L>0$ and consider the domains
\begin{align*}
	\dom&:=(-\infty,\infty)\times(0,\pi),\quad
	&\domm&:=(-\infty,0)\times(0,\pi),\quad
	&\domp&:=(0,\infty)\times(0,\pi),\\
	\domL&:=(-L,L)\times(0,\pi),\quad
	&\domLm&:=(-L,0)\times(0,\pi),\quad
	&\domLp&:=(0,L)\times(0,\pi).
\end{align*}
On $H^1_0(D)$, $D=\dom,\domL$ we work with the scalar product $\langle u,u\tf\rangle_{H^1_0(D)}:=\langle\nabla u,\nabla u\tf\rangle_{D}$.
For the bounded domain $\domL$ the equivalence of $\langle\cdot,\cdot\rangle_{H^1_0(\domL}$ to the standard $H^1(\domL)$-scalar product is well known.
For the unbounded domain $\dom$ this equivalence requires a short discussion:
Let
\begin{align*}
	\bfy_m(y):=\sqrt{\frac{2}{\pi}}\sin(m y), \quad m\in\setN.
\end{align*}
and recall that each $u\in H^1_0(\dom)$ and $u\in H^1_0(\domL)$ admits a representation
\begin{align}\label{eq:fourier}
	u(x,y)=\sum_{m\in\setN} u_m(x)\bfy_m(y), \quad
	u_m(x):=\langle u(x,\cdot),\bfy_m\rangle_{(0,\pi)}
\end{align}
with
\begin{align*}
\|u\|_{H^1(\dom)}^2 &= \sum_{m\in\setN} \|\partial_x u_m\|_{L^2(\setR)}^2
	+(\ev_m^2+1)\|u_m\|_{L^2(\setR)}^2
	\qquad\text{and}\\
\|u\|_{H^1(\domL)}^2 &= \sum_{m\in\setN} \|\partial_x u_m\|_{L^2(-L,L)}^2
	+(\ev_m^2+1)\|u_m\|_{L^2(-L,L)}^2
\end{align*}
respectively, where
\begin{align*}
	\ev_m:=m, \quad m\in\setN.
\end{align*}
It follows that $\|u\|_{H^1_0(\dom)}^2 \geq \frac{1}{2} \|u\|_{H^1(\dom)}^2$.
As usual, we consider any subspaces of $H^1_0(\dom)$ and $H^1_0(\domL)$ to be equipped with their inherited scalar product.

Let \(\sigma\) be constant on $\dompm$ with values \(\sigma_{-}:=\sigma|_{\domm}<0\) and \(\sigma_{+}:=\sigma|_{\domp}>0\).
Let $f\in L^2(\domL)$ and identify $f$ with its continuation by zero to $\dom$.
We consider the following two model problems:
\begin{subequations}
\begin{align}
\label{al:org_prob}
\begin{aligned}
	\text{Find } u \in H^1_0(\dom) \text{ such that }
	-\div(\sigma\nabla u)=f \text{ in }\dom; 
\end{aligned}\\
\label{al:org_prob_dir}
\begin{aligned}
	\text{Find } u \in H^1_0(\domL) \text{ such that }
	-\div(\sigma\nabla u)=f \text{ in }\domL,
\end{aligned}
\end{align}
\end{subequations}
and their variational formulations:
\begin{subequations}
\begin{align}
\label{al:org_prob_var}
	\text{Find } u \in H^1_0(\dom) \text{ such that }
	\ses_{\dom}(u,u\tf)
	&=\langle f,u\tf\rangle_\dom
	\text{ for all }u\tf\in H^1_0(\dom);\\
\label{al:org_prob_dir_var}
	\text{Find } u \in H^1_0(\domL) \text{ such that }
	\ses_{\domL}(u,u\tf)
	&=\langle f,u\tf\rangle_{\domL}
	\text{ for all }u\tf\in H^1_0(\domL),
\end{align}
\end{subequations}
with the corresponding sesquilinear forms
\begin{align*}
	\ses_D(u,u\tf):=\langle \sigma\nabla u,\nabla u\tf\rangle_{D},
	\quad D=\dom,\domL.
\end{align*}
Furthermore, let $\op\in\Lspace(H^1_0(\dom))$, $\opL\in\Lspace(H^1_0(\domL))$ be the associated operators defined by
\begin{subequations}
\begin{align}
	\langle \op u,u\tf \rangle_{H^1_0(\dom)}&=\ses_{\dom}(u,u\tf) \text{ for all }u,u\tf\in H^1_0(\dom),\\
	\langle \opL u,u\tf \rangle_{H^1_0(\domL)}&=\ses_{\domL}(u,u\tf) \text{ for all }u,u\tf\in H^1_0(\domL).
\end{align}
\end{subequations}
To specify the approxmiations of the former problems
let $P_1$ be the space of polynomials in one variable of order lower equal than one.
Thence let
\begin{align*}
	\Ny&\in\setN,\quad
	\hy:=\pi/\Ny,\quad
	\gridy_m:=\hy m \,\text{ for }\, m=1,\dots,\Ny,
\end{align*}
and
\begin{align*}
	\Vy&:=\{u\in H^1_0(0,\pi)\colon u|_{(\gridy_m,\gridy_{m+1})}\in P_1
	\text{ for all }
	m=0,\dots,\Ny-1\}.
\end{align*}
To discretize \eqref{al:org_prob_var} consider
\begin{align*}
	\hxpm&>0,\quad
	\gridx_n:=\hxp n, \,\text{ for }\, n=0,1,\dots;\quad
	\gridx_n:=\hxm n, \,\text{ for }\, n=-1,-2,\dots,
\end{align*}
and
\begin{align*}
	\Vx&:=\{u\in H^1(\setR)\colon
	u|_{(\gridx_n,\gridx_{n+1})}\in P_1 \text{ for all } n=\dots,-1,0,1,\dots\},
\end{align*}
and to discretize \eqref{al:org_prob_dir_var} let
\begin{align*}
	\Nxpm\in\setN,\quad
	\hxpm:=L/\Nxpm,\quad
	\gridxL_n:=\hxp n, \,\text{ for }\, n=0,\dots,\Nxp;\quad
	\gridxL_n:=\hxm n, \,\text{ for }\, n=-1,\dots,-\Nxm,
\end{align*}
and
\begin{align*}
	\VxL&:=\{u\in H^1_0(-L,L)\colon
	u|_{(\gridxL_n,\gridxL_{n+1})}\in P_1 \text{ for all } n=-\Nxm,\dots,\Nxp-1\}.
\end{align*}
Note that $\Vx$ is \emph{one} space which depends on \emph{both} parameters $\hxp$ and $\hxm$.
The same applies to $\VxL$.
Consequently we consider the Galerkin approximations of \eqref{al:org_prob_var} and \eqref{al:org_prob_dir_var} with discrete tensor product spaces
$\Vfull\subset H^1_0(\dom)$, $\VfullL\subset H^1_0(\domL)$:
\begin{subequations}
\begin{align}
\label{al:org_prob_dis}
\begin{aligned}
	\text{Find } u \in \Vfull \text{ such that }
	\ses_{\dom}(u,u\tf)
	=\langle f,u\tf\rangle_\dom
	\text{ for all }u\tf\in\Vfull;
\end{aligned}\\
\label{al:org_prob_dir_dis}
\begin{aligned}
	\text{Find } u \in \VfullL \text{ such that }
	\ses_{\domL}(u,u\tf)
	=\langle f,u\tf\rangle_{\domL}
	\text{ for all }u\tf\in\VfullL.
\end{aligned}
\end{align}
\end{subequations}
Let $\op_{\hxpm,\hy}\in\Lspace(\Vfull)$ and $\opL_{\hxpm,\hy}\in\Lspace(\VfullL)$ be the associated operators defined by
\begin{subequations}
\begin{align}
	\langle \op_{\hxpm,\hy} u,u\tf \rangle_{H^1_0(\dom)}&=\ses_{\dom}(u,u\tf) \text{ for all }u,u\tf\in\Vfull,\\
	\langle \opL_{\hxpm,\hy} u,u\tf \rangle_{H^1_0(\domL)}&=\ses_{\domL}(u,u\tf) \text{ for all }u,u\tf\in\VfullL.
\end{align}
\end{subequations}
We will see that the constrast $\kappa$ (of $\sigma)$ and the ratios of the meshes sizes
\begin{align*}
	\contrast:=\frac{\sigma_{+}}{\sigma_{-}}, \qquad 
	\ratio:=\frac{\hxp}{\hxm}, \qquad 
	\ratioy:=\frac{\hy}{\hxm},
\end{align*}
will play a crucial role in the stability analysis.
Since problem~\eqref{al:org_prob_var} is posed on an unbounded domain its discretization~\eqref{al:org_prob_dis} is rather theoretical, but it allows us to perform a very explicit analysis.
On the other hand, problem~\eqref{al:org_prob_dir_var} is posed on a bounded domain and hence its discretization~\eqref{al:org_prob_dir_dis} is computationally feasable, but its analyis is a bit more technical.
Indeed, the second setting \eqref{al:org_prob_dir_var}/\eqref{al:org_prob_dir_dis} can be considered as an approximation of \eqref{al:org_prob_var}/\eqref{al:org_prob_dis} by a truncation of the domain $\dom$ to $\domL$.
During the course of our analysis we will repeatedly use tensor product functions for which we apply the following notation in general:
\begin{align*}
	\uten(x,y)=\ux(x)\uy(y).
\end{align*}
In addition, let $\hfx_n\in\Vx,\hfxL_n\in\VxL,\hfy_m\in\Vy$ be the nodal basis functions defined by
\begin{align*}
	\hfx_n(\gridx_l)&=\delta_{nl},\quad n,l\in\setZ,\\
	\hfxL_n(\gridxL_l)&=\delta_{nl},\quad n,l=-\Nxm+1,\dots,\Nxp-1,\\
	\hfy_m(\gridy_l)&=\delta_{ml},\quad m,l=1,\dots,\Ny-1.
\end{align*}

\section{Stability analysis}\label{sec:discretization}
In this section we investigate the discretizations of \eqref{al:org_prob_var} and \eqref{al:org_prob_dir_var}.

\subsection{Unbounded domain}

To study the discretization of \eqref{al:org_prob_var} we first analyze its well-posedness and subsequently analyze a semi-discretization before treating the full discretization.

\subsubsection{Well-posedness analysis}\label{subsubsec:cont}

We start by discussing the well-posed of \eqref{al:org_prob} to ensure that we have chosen a meaningful problem.
Secondly, our analysis will serve as recipe for the forthcoming analysis dealing with the discretizations of \eqref{al:org_prob}.
Let
\begin{align*}
	X_{\pm}:=\{u\in H^1_0(\dom)\colon u|_{\dom_{\mp}}=0\}\text{ and } X_{0}:=(X_- \oplus X_+)^{\bot}.
\end{align*}
\begin{lemma}\label{lem:decomposition_cont}
The space $H^1_0(\dom)$ admits an orthogonal decomposition
\begin{align*}
	H^1_0(\dom)=X_- \oplus^\bot X_0 \oplus^\bot X_+
\end{align*}
where $X_0$ is spanned by the orthonormal basis
$\big( \frac{1}{\sqrt{2\ev_m}}
e^{-\ev_m|x|} \otimes\bfy_m(y)\big)_{m\in\setN}$.
\end{lemma}
\begin{proof}
Let $u\in X_0$.
By means of \eqref{eq:fourier} we can write $u=\sum_{m\in\setN} u_m\otimes\bfy_m$ and it follows that $u_m$ solves $-\partial_x\partial_x u_m+m^2 u_m=0$ in $\setR^\pm$.
Thus $u_m|_{\setR^+}(x)=c_1^+ e^{-\lambda_m x}+c_2^+ e^{\lambda_m x}$ and
$u_m|_{\setR^-}(x)=c_1^- e^{\lambda_m x}+c_2^- e^{-\lambda_m x}$ with constants $c_1^\pm,c_2^\pm\in\setR$.
Since $u_m\in H^1(\setR)$ it follows that $c_2^+=c_2^-=0$ and the continuity at the origin demands $c_1^+=c_1^-=:c$, i.e., $u_m(x)=c e^{-\lambda_m |x|}$.
The equality $c=1/(4\lambda_m)$ follows from a simple computation.
\end{proof}
\begin{lemma}\label{lem:blockoperator_cont}
The operator $\op$ is block diagonal with respect to the decomposition of \Cref{lem:decomposition_cont}.
The blocks corresponding to $X_-$, $X_+$ and $X_0$ equal the identity times $\sigma_-$, $\sigma_+$ and $\frac{\sigma_++\sigma_-}{2}=\sigma_-\frac{1+\contrast}{2}$ respectively.
Non-verbally: For each
$u_-,u_-\tf\in X_-$,
$u_0,u_0\tf\in X_0$,
$u_+,u_+\tf\in X_+$
it holds that
\begin{align*}
\ses_{\dom}(u_-+u_0+u_+,u_-\tf+u_0\tf+u_+\tf)
&=\sigma_-\langle  u_-,u_-\tf \rangle_{H^1_0(\dom)}
+\sigma_-\frac{1+\contrast}{2} \langle  u_0,u_0\tf \rangle_{H^1_0(\dom)}
+\sigma_+\langle  u_+,u_+\tf \rangle_{H^1_0(\dom)}.
\end{align*}
\end{lemma}
\begin{proof}
Since the supports of functions in $X_-$ and $X_+$ are disjoint it holds that $\ses_{\dom}(u_+,u_-)=\ses_{\dom}(u_-,u_+)=0$ for all $u_\pm\in X_\pm$.
We note that
\begin{align*}
\ses_{\dom}(u_0,u_\pm)
=\langle \sigma \nabla u_0, \nabla u_\pm\rangle_{\dom}
=\langle \sigma_\pm \nabla u_0, \nabla u_\pm\rangle_{\dompm}
=\sigma_\pm \langle \nabla u_0, \nabla u_\pm\rangle_{\dom}
=\sigma_\pm \langle u_0, u_\pm \rangle_{H^1_0(\dom)}
=0
\end{align*}
and likewise $\ses_{\dom}(u_\pm,u_0)=0$ for all $u_0\in X_0,u_\pm\in X_\pm$.
It remains to prove $\ses_{\dom}(u_0,u_0\tf)=\sigma_-\frac{1+\contrast}{2} \langle  u_0,u_0\tf \rangle_{H^1_0(\dom)}$, which follows from plugging in the orthonormal basis functions given in \Cref{lem:decomposition_cont}.
\end{proof}
\begin{corollary}\label{cor:wellposed}
If $\contrast\neq-1$, then $\op$ is bijective and $\|\op^{-1}\|_{\Lspace(H^1_0(\dom))}\leq \frac{1}{|\sigma_-|\min\{1,\contrast,\left|\frac{1+\contrast}{2}\right|\}}$.
\end{corollary}
\begin{proof}
From \Cref{lem:blockoperator_cont} we know that on each of the three subspaces \(\op\) is the identity operator multiplied by the constants \(\sigma_{-}\), \(\sigma_{+}\) and \(\sigma_{-}\frac{1+\contrast}{2}\) respectively. This implies that these are the only eigenvalues of the operator. It is therefore bijective as long as all of the three constant are nonzero, which is the case if \(\contrast\neq-1\), and the norm of the inverse is the reciprocal of the smallest of the three eigenvalues.
\end{proof}

\subsubsection{Semi discretization}\label{subsubsec:semi}

In this section we consider a semi discretization of problem \eqref{al:org_prob} by means of Galerkin spaces \(\Vsemi\):
\begin{align}\label{al:semidisc_prob}
	\text{Find } u \in \Vsemi \text{ such that }
	\ses_{\dom}(u,u\tf)=\langle f,u\tf\rangle_{\dom}
	\text{ for all } u\tf \in \Vsemi.
\end{align}
Let $\op_{\hxpm}\in\Lspace(\Vsemi)$ be the associated operator defined by
\begin{align}
	\langle \op_{\hxpm} u,u\tf \rangle_{H^1_0(\dom)}&=\ses_{\dom}(u,u\tf) \text{ for all }u,u\tf\in\Vsemi.
\end{align}
Note that by means of \eqref{eq:fourier} we can express:
\begin{align}\label{eq:ses_fourier_ortho}
	\ses_{\dom}(u,u\tf) &= \sum_{m\in\setN} \langle \sigma \partial_x u_m, \partial_x u_m\tf \rangle_{\setR}+\ev_m^2\langle \sigma u_m, u_m\tf \rangle_{\setR}.
\end{align}
In addition, for $u_m=\sum_{n\in\setZ} \alpha_n \hfx_n$, $u_m\tf=\sum_{n\in\setZ} \alpha_n\tf \hfx_n$, $\alpha_n,\alpha_n\tf\in\setR$ it follows that
\begin{align*}
	\langle \sigma \partial_x u_m, \partial_x u_m\tf \rangle_{\setR}
	+\ev_m^2\langle \sigma u_m, u_m\tf \rangle_{\setR}
	= \sum_{n,n'\in\setZ} \alpha_{n'}\tf \matAm_{n',n}  \alpha_n,
\end{align*}
where
\begin{align*}
\matAm:=\left(
    \begin{array}{ccccc}
         \ddots&\ddots&&&  \\
		 \ddots&2\sigma_{-}\amm&\sigma_{-}\bmm&&  \\
		 &\sigma_{-}\bmm&\sigma_{-}\amm+\sigma_{+}\amp&\sigma_{+}\bmp&  \\
         &&\sigma_{+}\bmp&2\sigma_{+}\amp&\ddots  \\
         &&&\ddots&\ddots
    \end{array}\right)
\end{align*}
with
\begin{subequations}\label{eq:ab}
\begin{align}
	\ampm&:=\phantom{-}\frac{1}{\hxpm}+\ev_m^{2}\frac{1}{3}\hxpm,\\
	\bmpm&:=-\frac{1}{\hxpm}+\ev_m^{2}\frac{1}{6}\hxpm.
\end{align}
\end{subequations}
Note that since $(\hfx_n)_{n\in\setZ}$ is not a Hilbert space basis the former expansion is only justified under certain decay conditions on $(\alpha_n)_{n\in\setZ}, (\alpha_n\tf)_{n\in\setZ}$.
However, this will not pose any problem for our forthcoming analysis.

We note that the analysis of the case $\bmpm=0$ is rather trivial, because thence $\matAm$ is diagonal.
Thus to unify our formulas we introduce the following case-wise definition
\begin{align*}
	\mumpmI&:=\begin{cases}
		\frac{1}{\bmpm}
		\left(-\ampm+\sqrt{(\ampm)^2-(\bmpm)^2}\right),&\bmpm\neq0,\\
		0,&\bmpm=0,
			\end{cases}\\
	\mumpmII&:=\begin{cases}
		\frac{1}{\bmpm}
		\left(-\ampm-\sqrt{(\ampm)^2-(\bmpm)^2}\right),&\bmpm\neq0,\\
		1,&\bmpm=0.
			\end{cases}
			\end{align*}
Henceforth we will only discuss the case $\bmpm\neq0$ and just note that the statements of all Lemmas and Theorems also hold for the case $\bmpm=0$.
Indeed $\mumpmI,\mumpmII$ are the roots of the polynomial $\mu\mapsto \bmpm\mu^2+2\ampm\mu+\bmpm$.

Note that it follows from $\mumpmI\mumpmII=1$ and \((\ampm)^2-(\bmpm)^2>0\) that $|\mumpmI|<|\mumpmII|$ and therefore
\begin{align}
	\label{eq:mumpmIII}
	|\mumpmI|<1 \quad\text{and}\quad |\mumpmII|>1.
\end{align}
We introduce the abbreviation
\begin{align}
	\label{eq:mumpm}
	\mumpm:=\mumpmI
\end{align}

and note that per definition $|\mumpm|<1$.
As the next step we further exploit the former representations in the following two lemmas which are in analogy to \Cref{subsubsec:cont}.
To avoid misconceptions we emphasize that the $n$ in
$\mu_{m,\pm}^{\pm n}$
appearing \Cref{lem:decomposition_semi} is an actual power and not an index.
\begin{lemma}
\label{lem:decomposition_semi}
The space $\Vsemi$ admits an orthogonal decomposition
\begin{align*}
\Vsemi
=\Big(\cl{\Vx^-\otimes H^1_0(0,\pi)}\Big)
\oplus^\bot \Vx^0
\oplus^\bot \Big(\cl{\Vx^+\otimes H^1_0(0,\pi)}\Big),
\end{align*}
where $\Vx^-:=\{v\in\Vx\colon v|_{\setR_+}=0\}$ and
$\Vx^+:=\{v\in\Vx\colon v|_{\setR_-}=0\}$.
The subspace $\Vx^0$ is spanned by the orthonormal basis
$\big(\vm(x)\otimes \bfy_m(y)\big)_{m\in\setN}$, where
\begin{align*}
\vm(x)&:=\frac{1}{\sqrt{\bmm\mumm+\amm+\amp+\bmp\mump}}
\Big(\hfx_0(x)+\sum_{n\in\setN} \mump^n\hfx_n(x)
+\sum_{n\in\setN} \mumm^{n}\hfx_{-n}(x)\Big).
\end{align*}
\end{lemma}
\begin{proof}
Let \(v\in \Big(\cl{\Vx^-\otimes H^1_0(0,\pi)}
\oplus  \cl{\Vx^+\otimes H^1_0(0,\pi)}\Big)^{\bot}\).
By means of \eqref{eq:fourier} we can write $v=\sum_{m\in\setN} v_m\otimes\bfy_m$. We can write \(\vm=\sum_{n\in \setZ}\beta_{n}^{(m)}\phi_{n}\). The span of the functions \(\phi_{n}\otimes\bfy_{m}\) for \(m,n\in\setZ\), \(n \neq 0\) is dense in
\(\Big(\cl{\Vx^-\otimes H^1_0(0,\pi)}
\oplus  \cl{\Vx^+\otimes H^1_0(0,\pi)}\Big)\). By orthogonality we then have
\begin{align*}
    0&=\langle v,\phi_{\pm n}\otimes\bfy_{m}\rangle_{H^1_0(\dom)}=\langle \partial_x \vm,\partial_x\phi_{\pm n}\rangle_{\setR}+\ev_m^2\langle \vm,\phi_{\pm n}\rangle_{\setR}\\
    &=\langle \beta_{\pm n-1}^{(m)}\partial_x\phi_{\pm n-1},\partial_x\phi_{\pm n}\rangle_{\setR}+\langle \beta_{\pm n}^{(m)}\partial_x\phi_{\pm n},\partial_x\phi_{\pm n}\rangle_{\setR}+\langle \beta_{\pm n+1}^{(m)}\partial_x\phi_{\pm n+1},\partial_x\phi_{\pm n}\rangle_{\setR}\\
    &+\ev_{m}^{2}\left(\langle \beta_{\pm n-1}^{(m)}\phi_{\pm n-1},\phi_{\pm n}\rangle_{\setR}+\langle \beta_{\pm n}^{(m)}\phi_{\pm n},\phi_{\pm n}\rangle_{\setR}+\langle \beta_{\pm n+1}^{(m)}\phi_{\pm n+1},\phi_{\pm n}\rangle_{\setR}\right)\\
	&=\sigma_\pm\big(\beta_{\pm n-1}^{(m)}\bmpm+2\beta_{\pm n}^{(m)}\ampm+\beta_{\pm n+1}^{(m)}\bmpm\big).
\end{align*}
Solving this three-term recurrence relation and recalling \eqref{eq:mumpm}, \eqref{eq:mumpmIII} we obtain that
\begin{align*}
    \beta_{n}^{(m)}=\beta_{0}^{(m)}\mump^{n} \text{ and } \beta_{-n}^{(m)}=\beta_{0}^{(m)}\mumm^{n} \quad\forall n\in\setN.
\end{align*}
Finally we compute the normalization constant by
\begin{align*}
    1&=\langle \partial_x \vm, \partial_x \vm\rangle_{\setR}+\ev_{m}^{2}\langle\vm, \vm\rangle_{\setR}
	=\sum_{n\in\setZ} \beta^{(m)}_{n}\big(\langle\partial_x\phi_{n}, \partial_x \vm\rangle_{\setR}+\ev_{m}^{2}\langle\phi_{n},\vm\rangle_{\setR}\big)\\
	&=\beta^{(m)}_{0}\big(\langle\partial_x\phi_{0}, \partial_x \vm\rangle_{\setR}+\ev_{m}^{2}\langle\phi_{0},\vm\rangle_{\setR}\big)\\
	&=(\beta^{(m)}_{0})^2\left(\mumm \bmm+\amm+\amp+\mump\bmp\right),
\end{align*}
i.e.,
\begin{align*}
    \beta^{(m)}_{0}&=\frac{1}{\sqrt{\mumm \bmm+\amm+\amp+\mump\bmp}}.
\end{align*}
Note that this calculation also ensures that \(\vm\otimes \bfy_{m}\) has finite \(H^1_0(\dom)\)-norm, i.e., \(\vm\otimes \bfy_{m}\in H^1_0(\dom)\) is well defined.
\end{proof}
\begin{lemma}
\label{lem:blockoperator_semi}
The operator $\op_{\hxpm}$ is block diagonal with respect to the orthogonal decomposition given in Lemma \ref{lem:decomposition_semi}.
The blocks corresponding to $\cl{\Vx^-\otimes H^1_0(0,\pi)}$ and $\cl{\Vx^+\otimes H^1_0(0,\pi)}$ equal the identity times $\sigma_-$ and $\sigma_+$ respectively.
The block corresponding to $\Vx^0$ is diagonal with respect to the basis given in Lemma \ref{lem:decomposition_semi} and the diagonal entries are given by
\begin{align}\label{eq:dm}
\dm:=
\frac{\sigma_-\bmm\mumm+\sigma_-\amm+\sigma_+\amp+\sigma_+\bmp\mump}%
{\phantom{\sigma_-}\bmm\mumm+\phantom{\sigma_-}\amm+\phantom{\sigma_+}\amp+\phantom{\sigma_+}\bmp\mump},
\quad m\in\setN.
\end{align}
Non-verbally: For each
$u_-,u_-\tf\in\cl{\Vx^-\otimes H^1_0(0,\pi)}$,
$u_0,u_0\tf\in\Vx^0$,
$u_+,u_+\tf\in\cl{\Vx^+\otimes H^1_0(0,\pi)}$
and $u_0=\sum_{m\in\setN} \beta_m \vm\otimes\bfy_m$,
$u_0\tf=\sum_{m\in\setN} \beta_m\tf \vm\otimes\bfy_m$,
$(\beta_m)_{m\in\setN},(\beta_m\tf)_{m\in\setN}\in\ell^2(\setN)$
it holds that
\begin{align*}
\ses_{\dom}(u_-+u_0+u_+,u_-\tf+u_0\tf+u_+\tf)
&=\sigma_-\langle  u_-,u_-\tf \rangle_{H^1_0(\dom)}
+\sum_{m\in\setN} \dm \beta_m\beta_m\tf
+\sigma_+\langle  u_+,u_+\tf \rangle_{H^1_0(\dom)}.
\end{align*}
\end{lemma}
\begin{proof}
The arguments are essentially the same as in the continuous case (see \Cref{lem:blockoperator_cont}) where we now use \Cref{lem:decomposition_semi} instead of \Cref{lem:decomposition_cont}. With this we directly get \(\ses_{\dom}(u_{-},u_{+}\tf)=\ses_{\dom}(u_{+},u_{-}\tf)=0\) and that \(\op_{\hxpm}\) is the identity times \(\sigma_{+}\) and \(\sigma_{-}\) on $\cl{\Vx^-\otimes H^1_0(0,\pi)}$ and $\cl{\Vx^+\otimes H^1_0(0,\pi)}$ respectively. By the orthogonality of the decomposition we also get \(\ses_{\dom}(u_{0},u_{\pm}\tf)=\sigma_{\pm}\langle u_{0},u_{\pm}\tf\rangle_{H_{0}^{1}(\dom)}=0\) as in the continuous case. It remains to show, that operator is diagonal on \(\Vx^0\) with the claimed values. That it is indeed diagonal follows directly from  the decomposition in \eqref{eq:ses_fourier_ortho}. To compute the values we use the same calculation as for the normalization constant and get
\begin{align*}
    \ses_{\dom}(\vm\otimes\bfy_{m},\vm\otimes\bfy_{m})&=\frac{\sum_{n\in \setN}\ses_{\dom}(\mumm^{n}\phi_{-n}\otimes\bfy_{m},\vm\otimes\bfy_{m})+\sum_{n\in \setN}\ses_{\dom}(\mump^{n}\phi_{n}\otimes\bfy_{m},\vm\otimes\bfy_{m})}{\bmm\mumm+\amm+\amp+\bmp\mump}\\
    &+\frac{\ses_{\dom}(\phi_{0}\otimes\bfy_{m},\vm\otimes\bfy_{m})}{\bmm\mumm+\amm+\amp+\bmp\mump}\\
 &=\frac{\ses_{\dom}(\phi_{0}\otimes\bfy_{m},\vm\otimes\bfy_{m})}{\bmm\mumm+\amm+\amp+\bmp\mump}\\
    &=\frac{\sigma_-\bmm\mumm+\sigma_-\amm+\sigma_+\amp+\sigma_+\bmp\mump}%
{\phantom{\sigma_-}\bmm\mumm+\phantom{\sigma_-}\amm+\phantom{\sigma_+}\amp+\phantom{\sigma_+}\bmp\mump},
\end{align*}
where we exploited the orthogonality properties of $\vm$.
\end{proof}
We observe that in contrast to \Cref{lem:blockoperator_cont} the block corresponding to $\Vx^0$ is not a multiple of the identity, but still diagonal.
To analyze the diagonal entries $\dm$ we introduce the function
\begin{align*}
\mf_{\contrast,\ratio}(t):=
\frac{1+\frac{\contrast\sqrt{\ratio^2t^2+12})}{\sqrt{t^2+12}}}%
{1+\frac{\sqrt{\ratio^2t^2+12})}{\sqrt{t^2+12}}}.
\end{align*}
\begin{lemma}\label{lem:diagonal}
The diagonal entries $\dm$ defined in \eqref{eq:dm} satisfy 
$\dm=\sigma_-\mf_{\contrast,\ratio}(\ev_m\hxm)$.
\end{lemma}
\begin{proof}
To start with, plugging in the definitions \eqref{eq:mumpm}, \eqref{eq:ab} of $\mumpm$ and $\bmpm, \ampm$ respectively yields that
\begin{align*}
    \bmpm\mumpm+\ampm&=\sqrt{\left(\ampm\right)^{2}-\left(\bmpm\right)^{2}}=\ev_m\sqrt{1+\frac{1}{12}\ev_m^2\hxpm^{2}}.
\end{align*}
Inserting this into the definition \eqref{eq:dm} of \(\dm\) we obtain that
\begin{align*}
    \dm&=\frac{\sigma_-(\bmm\mumm+\amm)+\sigma_+(\amp+\bmp\mump)}%
{\phantom{\sigma_-}(\bmm\mumm+\amm)+\phantom{\sigma_+}(\amp+\bmp\mump)}=\frac{\sigma_-\sqrt{1+\frac{1}{12}\ev_m^2\hxm^{2}}+\sigma_+\sqrt{1+\frac{1}{12}\ev_m^2\hxp^{2}}}%
{\phantom{\sigma_-}\sqrt{1+\frac{1}{12}\ev_m^2\hxm^{2}}+\phantom{\sigma_+}\sqrt{1+\frac{1}{12}\ev_m^2\hxp^{2}}}\\
&=\sigma_-\frac{1+\contrast\frac{\sqrt{12+\ev_m^2\hxp^{2}}}{\sqrt{12+\ev_m^2\hxm^{2}}}}%
{1+\frac{\sqrt{12+\ev_m^2\hxp^{2}}}{\sqrt{12+\ev_m^2\hxm^{2}}}}=\sigma_-\frac{1+\contrast\frac{\sqrt{12+\ev_m^2\hxm^{2}\ratio^{2}}}{\sqrt{12+\ev_m^2\hxm^{2}}}}%
{1+\frac{\sqrt{12+\ev_m^2\hxm^{2}\ratio^{2}}}{\sqrt{12+\ev_m^2\hxm^{2}}}}=\sigma_-\mf_{\contrast,\ratio}(\ev_m\hxm),
\end{align*}
where we recall that $r=\hxp/\hxm$.
\end{proof}
In the following lemmas we analyze the function $\mf_{\contrast,\ratio}$.
\begin{lemma}\label{lem:rootexplicit}
If one of the following two conditions
\begin{align}
\label{eq:AssR}
\underbrace{|\contrast|<1 \text{ and } \ratio|\contrast|>1}_{(\ref{eq:AssR}a)}
\quad\text{ or }\quad
\underbrace{|\contrast|>1 \text{ and } \ratio|\contrast|<1}_{(\ref{eq:AssR}b)}
\end{align}
is satisfied, then the only root of $\mf_{\contrast,\ratio}$ in $[0,+\infty)$ is
\begin{align*}
\roott:=\sqrt{\frac{12(1-\contrast^2)}{\contrast^2\ratio^2-1}}.
\end{align*}
\end{lemma}
\begin{proof}
Since the denominator of $\mf_{\contrast,\ratio}$ ranges for $t\geq0$ between $2$ and $1+r$, it suffices to analyze its nominator. Because we only consider non-negative \(t\) and negative \(\contrast\) we get
\begin{align*}
    0=1+\frac{\contrast\sqrt{\ratio^2t^2+12})}{\sqrt{t^2+12}}
	\Leftrightarrow\sqrt{t^2+12}=|\kappa|\sqrt{\ratio^{2}t^2+12}
	&\Leftrightarrow t^2+12=\kappa^{2}(\ratio^{2}t^2+12)\\
	&\Leftrightarrow t^2(1-\kappa^{2}\ratio^{2})=12(\kappa^{2}-1).
\end{align*}
The condition \eqref{eq:AssR} now guarantees that \(\kappa^{2}-1\) and \(1-\kappa^{2}\ratio^{2}\) have the same sign so the only solution is
\(\roott=\sqrt{\frac{12(1-\contrast^2)}{\contrast^2\ratio^2-1}}\).
\end{proof}
\begin{lemma}\label{lem:lim_f_zero}
If \eqref{eq:AssR} is satisfied, then $\lim_{\hxm\to0+} \inf_{m\in\setN} |\mf_{\contrast,\ratio}(\ev_m\hxm)|=0$.
\end{lemma}
\begin{proof}
Write $\hxm=\frac{1}{l+\epsilon} \roott$ with $l\in\setN_0$ and $\epsilon\in[0,1)$.
Choose $m=l$, recall that $\ev_m=m$ and exploit that $\lim_{l\to+\infty}\frac{l}{l+\epsilon}\to1$ uniformly in $\epsilon\in[0,1)$.
Apply \Cref{lem:rootexplicit} and the continuity of $\mf_{\contrast,\ratio}$.
\end{proof}
\begin{lemma}\label{lem:liminf_f}
If one of the following two conditions
\begin{align}
\label{eq:AssR_vv}
\underbrace{|\contrast|<1 \text{ and } \ratio|\contrast|<1}_{(\ref{eq:AssR_vv}a)}
\quad\text{ or }\quad
\underbrace{|\contrast|>1 \text{ and } \ratio|\contrast|>1}_{(\ref{eq:AssR_vv}b)}
\end{align}
is satisfied, then
$\inf_{t\geq0} |\mf_{\contrast,\ratio}(t)|\geq\min\{\left|\frac{1+\kappa}{2}\right|,\left|\frac{1+r\kappa}{1+r}\right|\}$.
\end{lemma}
\begin{proof}
We note that
\(\mf_{\contrast,\ratio}'(t) = \frac{12(\contrast-1)(\ratio^2-1)t}{\sqrt{t^2+12}\sqrt{\ratio^2t^2+12}(\sqrt{\ratio^2t^2+12}+\sqrt{t^2+12})^2}\)
and hence \(\mf_{\contrast,\ratio}(t)\) is a monotone function.
This implies \[\mf_{\contrast,\ratio}(t)\in\left[\min\left\{\mf_{\contrast,\ratio}(0),\lim_{s\rightarrow+\infty}\mf_{\contrast,\ratio}(s)\right\},\max\left\{\mf_{\contrast,\ratio}(0),\lim_{s\rightarrow+\infty}\mf_{\contrast,\ratio}(s)\right\}\right]\quad \forall t\in [0,+\infty).\]
Computing these values we get \(\mf_{\contrast,\ratio}(t)\in[\min\{\frac{1+\kappa}{2},\frac{1+r\kappa}{1+r}\},\max\{\frac{1+\kappa}{2},\frac{1+r\kappa}{1+r}\}]\). The conditions now ensure that both values have the same sign so the absolute value of \(\mf_{\contrast,\ratio}(t)\) is always bigger than the minimal absolute value occurring in one of the bounds.
\end{proof}

Now we are in the position to conclude our analysis of the semi discretization \eqref{al:semidisc_prob} in the following theorem.
\begin{theorem}\label{thm:semi}
If \eqref{eq:AssR_vv}
is satisfied, then $\op_{\hxpm}^{-1}$ exists and satisfies
\begin{align*}
\|\op_{\hxpm}^{-1}\|_{\Lspace(\Vsemi)}\leq\frac{1}{|\sigma_-|\min\left\{1,|\contrast|,\left|\frac{1+\contrast}{2}\right|,\left|\frac{1+r\contrast}{1+r}\right|\right\}}
\end{align*}
for all $\hxm,\hxp=\ratio \hxm>0$.
Contrary, if \eqref{eq:AssR} is satisfied, then
\begin{align*}
\lim_{\substack{\hxm\to0+\\\hxp=\ratio\hxm}} \|\op_{\hxpm}^{-1}\|_{\Lspace(\Vsemi)}
=+\infty
\end{align*}
(where we define $\|\op_{\hxpm}^{-1}\|_{\Lspace(\Vsemi)}:=+\infty$, if $\op_{\hxpm}^{-1}$ does not exist),
and in particular: $\op_{\hxpm}$ admits a nontrivial kernel for each
\[
\hxm= \frac{1}{\ev_m}\roott
=\frac{1}{m}\sqrt{\frac{12(1-\contrast^{2})}{\contrast^{2}\ratio^{2}-1}}
,\quad
\hxp=\ratio\hxm,
\quad m\in \setN.
\]
\end{theorem}
\begin{proof}
By \Cref{lem:blockoperator_semi} we know that the eigenvalues of \(\op_{\hxpm}\) are \(\sigma_{+}\),\(\sigma_{-}\) and \(\dm\) for \(m \in \setN\). From \Cref{lem:lim_f_zero} we get a lower bound on the absolute values of \(\dm\) which implies the first statement of the theorem. Then from \Cref{lem:liminf_f} it follows that such a bound does not exist in the other case which implies the second statement and finally \Cref{lem:rootexplicit} gives us the precise values where we have a zero eigenvalue which implies a nontrivial kernel.
\end{proof}

\subsubsection{Full discretization}

Next we consider the full discretization \eqref{al:org_prob_dis} by means of the Galerkin spaces $\Vfull$.
To this end much of the analysis of \Cref{subsubsec:semi} can be repeated, but we have to replace the orthogonal basis $(\bfy_m)_{m\in\setN}$ of $H^1_0(0,\pi)$ by a suitable discrete orthogonal basis $(\bfydis_m)_{m=1,\dots,\Ny-1}$ of $\Vy$,
and as a general rule we denote respective modified quantities by the same symbol as previously used but with an additional hat.
Hence we consider the following eigenvalue problem:
\begin{align}\label{al:disc_eigv}
\text{Find } (\evtmp,w) \in \setR^+\times\Vy\setminus\{0\} \text{ such that }
\langle \partial_y w,\partial_y w\tf\rangle_{(0,2\pi)}
=\evtmp^2 \langle w,w\tf \rangle_{(0,2\pi)} \quad \forall w\tf \in \Vy.
\end{align}
To solve this problem we define
\begin{align*}
	(\matB)_{m',m}:=\langle \partial_y \hfy_m,\partial_y \hfy_{m'}\rangle_{(0,2\pi)}
	-\evtmp^2 \langle \hfy_m,\hfy_{m'} \rangle_{(0,2\pi)},
	\quad m,m'=1,\dots,\Ny-1.
\end{align*}
Then we use that \(\evtmp\) is an eigenvalue if and only if \(\matB\) has a zero eigenvalue.
It holds that
\begin{align*}
\matB&=\left(
\begin{array}{ccccc}
         2\adis&\bdis&&&  \\
         \bdis&2\adis&\ddots&&  \\
         &\ddots&\ddots&\ddots&  \\
         &&\ddots&2\adis&\bdis  \\
         &&&\bdis&2\adis
\end{array}\right),\qquad
\adis:= \frac{1}{\hy}-\evtmp^2\frac{1}{3} \hy,\qquad
\bdis:=-\frac{1}{\hy}-\evtmp^2\frac{1}{6} \hy.
\end{align*}
It follows that \cite{Boffi10,EkstroemSerraC18}
\begin{align*}
	\evtmp_{m}^2&=
	\frac{6}{\hy^{2}}\frac{1-\cos(\hy m)}{2+\cos(\hy m)},
	\quad m=1,\dots,\Ny-1
\end{align*}
with respective $L^2(0,\pi)$-normalized eigenfunctions
\begin{align*}
	\bfydis_m(y):=c\sum_{l=1}^{\Ny-1} \sin(m \hy l) \hfy_l(y),\quad
	c:=\Big\|\sum_{l=1}^{\Ny-1} \sin(m \hy l) \hfy_l\Big\|_{L^2(0,\pi)}^{-1},\quad
	m=1,\dots,\Ny-1.
\end{align*}
This leads us to introduce
\begin{align*}
\evdism:=\frac{\sqrt{6}}{\ratioy\hxm}\sqrt{\frac{1-\cos(\ratioy\hxm m)}{2+\cos(\ratioy\hxm m)}}
=\frac{\sqrt{6}}{\hy}\sqrt{\frac{1-\cos(\hy m)}{2+\cos(\hy m)}},
\quad m=1,\dots,\Ny-1.
\end{align*}
Hence we define respective modified quantities
\begin{align*}
\ampmhat&:=\frac{1}{\hxpm}+\evdism^{2}\frac{1}{3}\hxpm,\\
\bmpmhat&:=-\frac{1}{\hxpm}+\evdism^{2}\frac{1}{6}\hxpm,\\
\mumpmIhat&:=\frac{1}{\bmpmhat} \left(-\ampmhat+\sqrt{(\ampmhat)^2-(\bmpmhat)^2}\right),\\\mumpmIIhat&:=\frac{1}{\bmpmhat}	\left(-\ampmhat-\sqrt{(\ampmhat)^2-(\bmpmhat)^2}\right),\\
\mumpmhat&:=\begin{cases} \mumpmIhat,&\bmpmhat\neq0,\\ 0,&\bmpmhat=0, \end{cases},\\
\dmhat&:=\frac{\sigma_-\bmm\mummhat+\sigma_-\ammhat+\sigma_+\amphat+\sigma_+\bmphat\mumphat}%
{\phantom{\sigma_-}\bmmhat\mummhat+\phantom{\sigma_-}\ammhat+\phantom{\sigma_+}\amphat+\phantom{\sigma_+}\bmphat\mumphat},\\
\vmhat(x)&:=\frac{1}{\sqrt{\bmmhat\mummhat+\ammhat+\amphat+\bmphat\mumphat}}
\Big(\hfx_0(x)+\sum_{n\in\setZ^+} \mumphat^n\hfx_n(x)
+\sum_{n\in\setZ^-} \mummhat^{-n}\hfx_n(x)\Big)
\end{align*}
for $m=1,\dots,\Ny-1$.
The forthcoming two lemmas follow in analogy to \Cref{subsubsec:semi}.
\begin{lemma}\label{lem:decomposition_full}
The space $\Vfull$ admits an orthogonal decomposition
\begin{align*}
\Vfull
=\Big(\Vx^-\otimes\Vy\Big)
\oplus^\bot \Vxzhat
\oplus^\bot \Big(\Vx^+\otimes\Vy\Big),
\end{align*}
where $\Vx^-=\{v\in\Vx\colon v|_{\setR_+}=0\}$ and
$\Vx^+-=\{v\in\Vx\colon v|_{\setR_-}=0\}$.
The subspace $\Vxzhat$ is spanned by the orthonormal basis
$\big(\vmhat(x)\otimes \bfydis_m(y)\big)_{m=1,\dots,\Ny-1}$.
\end{lemma}
\begin{proof}
The proof can be obtained by following the steps of the proof of \Cref{lem:decomposition_semi} one-to-one.
\end{proof}
\begin{lemma}\label{lem:blockoperator_full}
The operator $\op_{\hxpm,\hy}$ is block diagonal with respect to the orthogonal decomposition given in \Cref{lem:decomposition_full}.
The blocks corresponding to $\Vx^-\otimes\Vy$ and $\Vx^+\otimes\Vy$ equal the identity times $\sigma_-$ and $\sigma_+$ respectively.
The block corresponding to $\Vxzhat$ is diagonal with respect to the basis given in \Cref{lem:decomposition_full} and the diagonal entries are given by $\dmhat$.
Non-verbally: For each
$u_-,u_-\tf\in\Vx^-\otimes\Vy$,
$u_0,u_0\tf\in\Vxzhat$,
$u_+,u_+\tf\in\Vx^+\otimes\Vy$
and $u_0=\sum_{m=1}^{\Ny-1} \beta_m \vmhat\otimes\bfydis_m$,
$u_0\tf=\sum_{m=1}^{\Ny-1} \beta_m\tf \vmhat\otimes\bfydis_m$,
$(\beta_m)_{m=1}^{\Ny-1},(\beta_m\tf)_{m=1}^{\Ny}\in\setR^{\Ny-1}$
it holds that
\begin{align*}
\ses_{\dom}(u_-+u_0+u_+,u_-\tf+u_0\tf+u_+\tf)
&=\sigma_-\langle  u_-,u_-\tf \rangle_{H^1_0(\dom)}
+\sum_{m=1}^{\Ny-1} \dmhat \beta_m\beta_m\tf
+\sigma_+\langle  u_+,u_+\tf \rangle_{H^1_0(\dom)}.
\end{align*}
\end{lemma}
\begin{proof}
The proof can be obtained by following the steps of the proof of \Cref{lem:blockoperator_semi} one-to-one. 
\end{proof}
To analyze the diagonal entries $\dmhat$ we define
\begin{align*}
	\mh_{\ratioy}(s):=\sqrt{\frac{6}{\ratioy^{2}}\frac{1-\cos(s)}{2+\cos(s)}}
\end{align*}
and introduce the following lemmas.
\begin{lemma}
It holds that 
$\evdism\hxm=\mh_{\ratioy}(\ratioy m\hxm)$
and
$\dmhat=\sigma_-\mf_{\contrast,\ratio}\big(\mh_{\ratioy}(\ratioy m\hxm)\big)$
for $m=1,\dots,\Ny-1$.
\end{lemma}
\begin{proof}
An elementary computation shows
\[\evdism^{2}\hxm^2=\frac{6}{\ratioy^{2}\hxm^{2}}\frac{1-\cos(\ratioy\hxm m)}{2+\cos(\ratioy\hxm m)}\hxm^{2}=\mh_{\ratioy}(\ratioy m\hxm)^2.\]
Then as in \Cref{lem:diagonal} it follows that $\dmhat=\sigma_-\mf_{\contrast,\ratio}\big(\evdism\hxm\big)$ and combining these two identities yields that
\(\dmhat=\sigma_-\mf_{\contrast,\ratio}\big(\mh_{\ratioy}(\ratioy m\hxm)\big).\)
\end{proof}
\begin{lemma}\label{lem:roots}
If one of the following two conditions
\begin{align}
\label{eq:AssRy}
\underbrace{|\contrast|<1 \text{ and } \ratio^2\contrast^2>1+\ratioy^{2}(1-\contrast^2)}_{(\ref{eq:AssRy}a)}
\quad\text{ or }\quad
\underbrace{|\contrast|>1 \text{ and } \ratio^2\contrast^2<1+\ratioy^{2}(1-\contrast^2)}_{(\ref{eq:AssRy}b)}
\end{align}
is satisfied, then the problem to find $s\in(0,\pi]$ such that $\mh_{\ratioy}(s)=\roott$ admits the unique solution
\begin{align*}
\roots:=\arccos\left(\frac{1-\frac{\ratioy^{2}\roott^{2}}{3}}{1+\frac{\ratioy^{2}\roott^{2}}{6}}\right)
=\arccos\left(1+\frac{6\ratioy^{2}(1-\contrast^{2})}{(1-\contrast^{2}\ratio^{2})-2\ratioy^{2}(1-\contrast^{2})}\right).
\end{align*}
\end{lemma}
\begin{proof}
First, one can check, that the inequalities guarantee that \(\frac{6\ratioy^{2}(1-\contrast^{2})}{(1-\contrast^{2}\ratio^{2})-2\ratioy^{2}(1-\contrast^{2})}\in (-2,0)\) so \(\roots\) is well defined. Now we compute
\begin{align*}
    \mh_{r_{y}}(\roots)^{2}=\frac{6}{\ratioy^{2}}\frac{1-\cos(\roots)}{2+\cos(\roots)}=\frac{6}{\ratioy^{2}}\frac{1-\frac{1-\frac{\ratioy^{2}\roott^{2}}{3}}{1+\frac{\ratioy^{2}\roott^{2}}{6}}}{2+\frac{1-\frac{\ratioy^{2}\roott^{2}}{3}}{1+\frac{\ratioy^{2}\roott^{2}}{6}}}=\frac{6}{\ratioy^{2}}\frac{3\frac{\ratioy^{2}\roott^{2}}{6}}{3}=\roott^{2}.
\end{align*}
\end{proof}
Before we formulate the next \Cref{lem:lim_fh_zero}, let us recall that $\hxm=\frac{\hy}{\ratioy}=\frac{\pi}{\ratioy\Ny}$.
\begin{lemma}\label{lem:lim_fh_zero}
If \eqref{eq:AssRy} is satisfied, then $\lim_{\Ny\to+\infty} \inf_{m\in\{1,\dots,\Ny-1\}}
\left|\mf_{\contrast,\ratio}\big(\mh_{\ratioy}(\frac{m\pi}{\Ny})\big)\right|=0$.
\end{lemma}
\begin{proof}
Because the rational numbers are dense in the real numbers, for each \(\epsilon>0\)
there exists \(\Ny_{\epsilon}\in \setN\) such that for all $\Ny\in\setN,\Ny>\Ny_\epsilon$ there exists \(m\in\{1,..,\Ny-1\}\) such that \(\left|\frac{m}{\Ny}\pi-\roots\right|<\epsilon\).
The theorem now follows from the continuity of \(\mf_{\contrast,\ratio}\circ\mh_{\ratioy}\) and \(\mf_{\contrast,\ratio}(\mh_{\ratioy}(\roots))=0\).
\end{proof}
\begin{lemma}\label{lem:liminf_f_full}
If one of the following two conditions
\begin{align}
\label{eq:AssRy_vv}
\underbrace{|\contrast|<1 \text{ and } \ratio^2\contrast^2<1+\ratioy^2(1-\contrast^2)}_{(\ref{eq:AssRy_vv}a)}
\quad\text{ or }\quad
\underbrace{|\contrast|>1 \text{ and } \ratio^2\contrast^2>1+\ratioy^2(1-\contrast^2)}_{(\ref{eq:AssRy_vv}b)}
\end{align}
is satisfied, then
$\inf_{s\in\setR} |\mf_{\contrast,\ratio}\big(\mh_{\ratioy}(s)\big)|\geq
\min\Big\{\left|\frac{1+\contrast}{2}\right|,\left|\frac{\sqrt{1+\ratioy^{2}}+\contrast\sqrt{\ratio^{2}+\ratioy^{2}}}{\sqrt{1+\ratioy^{2}}+\sqrt{\ratio^{2}+\ratioy^{2}}}\right|\Big\}>0$.
\end{lemma}
\begin{proof}
The reasoning is the same as in the proof of \Cref{lem:liminf_f}. The only difference in this case is, that \(\mh_{\ratioy}(s)\in[0,\frac{\sqrt{12}}{\ratioy}]\) so we do not consider the limit at infinity and get
\[\mf_{\contrast,\ratio}\left(\frac{\sqrt{12}}{\ratioy}\right)=\frac{1+\contrast\sqrt{\frac{\ratio^{2}+\ratioy^{2}}{1+\ratioy^{2}}}}{1+\sqrt{\frac{\ratio^{2}+\ratioy^{2}}{1+\ratioy^{2}}}}=\frac{\sqrt{1+\ratioy^{2}}+\contrast\sqrt{\ratio^{2}+\ratioy^{2}}}{\sqrt{1+\ratioy^{2}}+\sqrt{\ratio^{2}+\ratioy^{2}}}\]
instead. Again, the condition \eqref{eq:AssRy_vv} ensures, that both bounds have the same sign, so we can safely take the minimum of their absolute values.
\end{proof}

Now we are in the position to conclude our analysis of the full discretization \eqref{al:semidisc_prob} in the following theorem.
\begin{theorem}\label{thm:full}
If \eqref{eq:AssRy_vv} is satisfied, then $\op_{\hxpm,\hy}^{-1}$ exists and satisfies
\begin{align*}
\|\op_{\hxpm,\hy}^{-1}\|_{\Lspace(\Vfull)}\leq
\frac{1}{|\sigma_-|\min\left\{1,|\contrast|,\left|\frac{1+\contrast}{2}\right|,\left|\frac{\sqrt{1+\ratioy^{2}}+\contrast\sqrt{\ratio^{2}+\ratioy^{2}}}{\sqrt{1+\ratioy^{2}}+\sqrt{\ratio^{2}+\ratioy^{2}}}\right|\right\}}
\end{align*}
Contrary, if \eqref{eq:AssRy} is satisfied, then
\begin{align*}
\lim_{h_{-}\rightarrow 0} \|\op_{\hxpm,\hy}^{-1}\|_{\Lspace(\Vfull)}=+\infty
\end{align*}
(where we define $\|\op_{\hxpm,\hy}^{-1}\|_{\Lspace(\Vfull)}:=+\infty$, if $\op_{\hxpm,\hy}^{-1}$ does not exist),
and in particular
$\op_{\hxpm,\hy}$ admits a nontrivial kernel for each
\[
\hxm
=\frac{1}{m}\frac{1}{\ratioy}\roots=\frac{1}{m}\frac{1}{\ratioy}\arccos\left(1+\frac{6\ratioy^{2}(1-\contrast^{2})}{(1-\contrast^{2}\ratio^{2})-2\ratioy^{2}(1-\contrast^{2})}\right),\quad m\in\setN.
\]
Note that to simultaneously satisfy \(\frac{\pi}{\hy}\in\setN\) we can choose a particular \(\ratioy\) or \(\ratio\) such that \(\roots=\pi l/k, l,k\in\setN\), which yields that \(\frac{\pi}{\hy}\in\setN\) for \(m\in k\setN\).
\end{theorem}
\begin{proof}
The proof can be obtained by following the steps of the proof of \Cref{thm:semi} one-to-one and replacing needed lemmas by the corresponding lemmas from this section.
\end{proof}

\subsection{Bounded domain}

Now we consider a bounded domain for which actual numerical computations are possible.
Since the well-posedness analysis follows along the lines of \Cref{subsubsec:cont} we suffice ourselves with stating following lemma without proof.

\begin{lemma}\label{lem:wellposed_L}
The space $H^1_0(\domL)$ admits an orthogonal decomposition
\(
	H^1_0(\domL)=\XL_- \oplus^\bot \XL_0 \oplus^\bot \XL_+
\)
where
\(
	\XL_{\pm}:=\{u\in H^1_0(\domL)\colon u|_{\domL_{\mp}}=0\}
\)
and $\XL_0$ is spanned by the orthonormal basis
\begin{align*}
	\left(\frac{1}{\sqrt{2\ev_{m}(e^{4\ev_{m}L}-1)}}
	\Big(e^{2\ev_m L}e^{-\ev_m|x|}-
	e^{\ev_m|x|} \Big) \otimes\bfy_m(y)\right)_{m\in\setN}.
\end{align*}
The operator \(\opL\) is block diagonal and the statements of \Cref{lem:blockoperator_cont} and \Cref{cor:wellposed} apply with \(\op,X_{\pm},X_{0}\) being replaced by \(\opL,\XL_{\pm},\XL_{0}\).
\end{lemma}

\subsubsection{Semi discretization}
We consider a semi discretization of problem \eqref{al:org_prob_dir} by means of Galerkin spaces \(\VsemiL\):
\begin{align}\label{al:semidisc_prob_L}
	\text{Find } u \in \VsemiL \text{ such that }
	\ses_{\domL}(u,u\tf)=\langle f,u\tf\rangle_{\domL}
	\text{ for all } u\tf \in \VsemiL.
\end{align}
Let $\opL_{\hxpm}\in\Lspace(\VsemiL)$ be the associated operator defined by
\begin{align}
	\langle \opL_{\hxpm} u,u\tf \rangle_{H^1_0(\domL)}&=\ses_{\domL}(u,u\tf) \text{ for all }u,u\tf\in\VsemiL.
\end{align}
As in the unbounded case we can use \eqref{eq:fourier} and derive a decomposition into Fourier modes as in \eqref{eq:ses_fourier_ortho}. Similarly we then also get that for $u_m=\sum_{n=-\Nxm}^{\Nxp} \alpha_n \hfx_n$, $u_m\tf=\sum_{n=-\Nxm}^{\Nxp} \alpha_n\tf \hfx_n$, $\alpha_n,\alpha_n\tf\in\setR$ it follows that
\begin{align*}
	\langle \sigma \partial_x u_m, \partial_x u_m\tf \rangle_{\setR}
	+\ev_m^2\langle \sigma u_m, u_m\tf \rangle_{\setR}
	= \sum_{n=-\Nxm}^{\Nxp}\sum_{n'=-\Nxm}^{\Nxp}\alpha_{n'}\tf \matAm_{n',n}  \alpha_n,
\end{align*}
where
\begin{align*}
\matAmL:=
\left(
    \begin{array}{ccccccc}
    2\sigma_{-}\amm&\sigma_{-}\bmm&&&&&  \\
         \sigma_{-}\bmm&\ddots&\ddots&&&&  \\
         &\ddots&2\sigma_{-}\amm&\sigma_{-}\bmm&&&  \\
         &&\sigma_{-}\bmm&\sigma_{-}\amm+\sigma_{+}\amp&\sigma_{+}\bmp&&  \\
         &&&\sigma_{+}\bmp&2\sigma_{+}\amp&\ddots&  \\
         &&&&\ddots&\ddots&\sigma_{+}\bmp\\
         &&&&&\sigma_{+}\bmp&2\sigma_{+}\amp
 	\end{array}\right).
\end{align*}
After these considerations we can again derive an orthogonal decomposition as follows.
\begin{lemma}\label{lem:decomposition_semi_L}
The space $\VsemiL$ admits an orthogonal decomposition
\begin{align*}
\VsemiL
=\Big(\VxL^-\otimes H^1_0(0,\pi)\Big)
\oplus^\bot \VxL^0
\oplus^\bot \Big(\VxL^+\otimes H^1_0(0,\pi)\Big),
\end{align*}
where $\VxL^-:=\{v\in\VxL\colon v|_{(0,L)}=0\}$ and
$\VxL^+-:=\{v\in\VxL\colon v|_{(-L,0)}=0\}$.
The subspace $\VxL^0$ is spanned by the orthonormal basis
$\big(\vmL(x)\otimes \bfy_m(y)\big)_{m\in\setN}$, where
\begin{align*}
\vmL(x)&:=\frac{1}{\sqrt{\frac{1+\nu_{m,-}^{\Nxm}}{1-\nu_{m,-}^{\Nxm}}\sqrt{(\amm)^2-(\bmm)^2}+\frac{1+\nu_{m,+}^{\Nxp}}{1-\nu_{m,+}^{\Nxp}}\sqrt{(\amp)^2-(\bmp)^2}}}
\Bigg(\hfx_0(x)\\
&\quad+\sum_{n=-\Nxm}^{-1}\frac{\mummI^{-n}-\nu_{m,-}^{\Nxm}\mummII^{- n}}{1-\nu_{m,-}^{\Nxm}}\hfx_n(x)
+\sum_{n=1}^{\Nxp}\frac{\mumpI^{n}-\nu_{m,+}^{\Nxp}\mumpII^{n}}{1-\nu_{m,+}^{\Nxp}}\hfx_n(x)\Bigg).
\end{align*}
Here \(\mumpmI, \mumpmII\) are as defined in \Cref{subsubsec:semi} and \(\nu_{m,\pm}:=\frac{\mumpmI}{\mumpmII}\).
\end{lemma}
\begin{proof}
The structure and beginning of the proof is the same as in \Cref{lem:decomposition_semi}. We again get the orthogonality of \(\VxL^-\otimes H^1_0(0,\pi)\) and \(\VxL^+\otimes H^1_0(0,\pi)\) because of non-intersecting supports and we get that the coefficients of \(\vmL\) in the finite element basis given by \(\vmL=\sum_{n=-\Nxm}^{\Nxp}\betamL_{n}\phi_{n}\) satisfy the system of equations
\begin{align*}
	\bmm \betamL_{n}+2\amm \betamL_{n+1}+\bmm \betamL_{n+2}&=0,\quad -\Nxm\leq n\leq -2,\\
	\bmp \betamL_{n}+2\amp \betamL_{n+1}+\bmp \betamL_{n+2}&=0,\quad 0\leq n\leq \Nxp-2.
\end{align*}
We now first treat the case when \(\bmm,\bmp\neq0\) and get that this system is solved by
\begin{align*}
	\betamL_{n}&=c_{1,-}\mummI^{-n}+c_{2,-}\mummII^{-n},\quad -\Nxm\leq n\leq 0,\\
	\betamL_{n}&=c_{1,+}\mumpI^{n}+c_{2,+}\mumpII^{n},\quad 0\leq n\leq \Nxp.
\end{align*}
where only the constants \(c_{1,\pm},c_{2,\pm}\) still have to be determined. To do this, we use that both equations hold for \(n=0\) and that \(\betamL_{-\Nxm}=\betamL_{\Nxp}=0\) because of the Dirichlet boundary conditions. This leads to the system
\begin{align*}
	c_{1,-}+c_{2,-}&=c_{1,+}+c_{2,+},\\
	c_{1,-}\mummI^{-\Nxm}+c_{2,-}\mummII^{-\Nxm}&=0,\\
	c_{1,+}\mumpI^{\Nxp}+c_{2,+}\mumpII^{\Nxp}&=0,
\end{align*}
with the solution 
\begin{align*}
    c_{1,+}=\frac{\betamL_{0}}{1-\nu_{m,+}^{\Nxp}},\quad c_{2,+}=-\frac{\betamL_{0}}{1-\nu_{m,+}^{\Nxp}}\nu_{m,+}^{\Nxp},\quad c_{1,-}=\frac{\betamL_{0}}{1-\nu_{m,-}^{\Nxm}},\quad c_{2,-}=-\frac{\betamL_{0}}{1-\nu_{m,-}^{\Nxm}}\nu_{m,+}^{\Nxm}.
\end{align*}
Here \(\betamL_{0}\) still has to be determined. This means, that 
\begin{align*}
	\betamL_{n}&=\betamL_{0}\frac{\mummI^{-n}-\nu_{m,-}^{\Nxm}\mummII^{- n}}{1-\nu_{m,-}^{\Nxm}},\quad -\Nxm\leq n\leq 0\\
	\betamL_{n}&=\betamL_{0}\frac{\mumpI^{n}-\nu_{m,+}^{\Nxp}\mumpII^{n}}{1-\nu_{m,+}^{\Nxp}},\quad 0\leq n\leq \Nxp.
\end{align*}
One can check that these formulas still hold in the case where \(\bmm\) or \(\bmp\) is zero if we then set \(\mumpmI\) to zero as we have done before and set \(\mumpmII\) to an arbitrary positive number. Finally we compute \(\betamL_{0}\) to normalize the solution. For this we only have to consider the scalar product with \(\phi_{0}\) because of orthogonality. This leads to
\begin{align*}
    1&=\langle \partial_x\vmL, \partial_x \vmL\rangle_{\setR}+\ev_{m}^{2}\langle\vmL,\vmL\rangle_{\setR}\\
	&=\sum_{n=\Nxm}^{\Nxp} \betamL_n \big(\langle \partial_x\phi_n, \partial_x \vmL\rangle_{\setR}+\ev_{m}^{2}\langle\phi_n,\vmL\rangle_{\setR}\big)\\
	&=\betamL_{0}\big(\langle \partial_x\phi_{0}, \partial_x \vmL\rangle_{\setR}+\ev_{m}^{2}\langle\phi_{0},\vmL\rangle_{\setR}\big)\\
	&=(\betamL_{0})^2\left(\frac{\mummI-\nu_{m,-}^{\Nxm}\mummII}{1-\nu_{m,-}^{\Nxm}}\bmm+\amm
    +\amp+\frac{\mumpI-\nu_{m,+}^{\Nxp}\mumpII}{1-\nu_{m,+}^{\Nxp}}\bmp\right)\\
    &=\left(\betamL_{0}\right)^{2}\left(\frac{(\mummI\bmm+\amm)-\nu_{m,-}^{\Nxm}(\mummII\bmm+\amm)}{1-\nu_{m,-}^{\Nxm}}\right.\\
	&\quad\phantom{\left(\betamL_{0}\right)^{2}}\left.+\frac{(\mumpI\bmp+\amp)-\nu_{m,+}^{\Nxp}(\mumpII\bmp+\amp)}{1-\nu_{m,+}^{\Nxp}}\right)\\
    &=\left(\betamL_{0}\right)^{2}\left(\frac{\sqrt{(\amm)^2-(\bmm)^2}+\nu_{m,-}^{\Nxm}\sqrt{(\amm)^2-(\bmm)^2}}{1-\nu_{m,-}^{\Nxm}}\right.\\
    &\quad\phantom{\left(\betamL_{0}\right)^{2}}\left.+\frac{\sqrt{(\amp)^2-(\bmp)^2}+\nu_{m,+}^{\Nxp}\sqrt{(\amp)^2-(\bmp)^2}}{1-\nu_{m,+}^{\Nxp}}\right)\\
	&=\left(\betamL_{0}\right)^{2}\left(\frac{1+\nu_{m,-}^{\Nxm}}{1-\nu_{m,-}^{\Nxm}}\sqrt{(\amm)^2-(\bmm)^2}+\frac{1+\nu_{m,+}^{\Nxp}}{1-\nu_{m,+}^{\Nxp}}\sqrt{(\amp)^2-(\bmp)^2}\right),
\end{align*}
i.e.,
\begin{align*}
    \betamL_{0}&=\frac{1}{\sqrt{\frac{1+\nu_{m,-}^{\Nxm}}{1-\nu_{m,-}^{\Nxm}}\sqrt{(\amm)^2-(\bmm)^2}+\frac{1+\nu_{m,+}^{\Nxp}}{1-\nu_{m,+}^{\Nxp}}\sqrt{(\amp)^2-(\bmp)^2}}}.
\end{align*}
\end{proof}
\begin{lemma}
\label{lem:blockoperator_semi_L}
The operator $\opL_{\hxpm}$ is block diagonal with respect to the orthogonal decomposition given in \Cref{lem:decomposition_semi_L}.
The blocks corresponding to $\VxL^-\otimes H^1_0(0,\pi)$ and $\VxL^+\otimes H^1_0(0,\pi)$ equal the identity times $\sigma_-$ and $\sigma_+$ respectively.
The block corresponding to $\VxL^0$ is diagonal with respect to the basis given in \Cref{lem:decomposition_semi_L} and the diagonal entries are given by
\begin{align}\label{eq:dm_L}
\dmL:=
\frac{\sigma_{-}\frac{1+\nu_{m,-}^{\Nxm}}{1-\nu_{m,-}^{\Nxm}}\sqrt{(\amm)^2-(\bmm)^2}+\sigma_{+}\frac{1+\nu_{m,+}^{\Nxp}}{1-\nu_{m,+}^{\Nxp}}\sqrt{(\amp)^2-(\bmp)^2}}%
{\phantom{\sigma_{-}}\frac{1+\nu_{m,-}^{\Nxm}}{1-\nu_{m,-}^{\Nxm}}\sqrt{(\amm)^2-(\bmm)^2}+\phantom{\sigma_{+}}\frac{1+\nu_{m,+}^{\Nxp}}{1-\nu_{m,+}^{\Nxp}}\sqrt{(\amp)^2-(\bmp)^2}},
\quad m\in\setN.
\end{align}
Non-verbally: For each
$u_-,u_-\tf\in\VxL^-\otimes H^1_0(0,\pi)$,
$u_0,u_0\tf\in\VxL^0$,
$u_+,u_+\tf\in\VxL^+\otimes H^1_0(0,\pi)$
and $u_0=\sum_{m\in\setN} \beta_m \vmL\otimes\bfy_m$,
$u_0\tf=\sum_{m\in\setN} \beta_m\tf \vmL\otimes\bfy_m$,
$(\beta_m)_{m\in\setN},(\beta_m\tf)_{m\in\setN}\in\ell^2(\setN)$
it holds that
\begin{align*}
\ses_{\domL}(u_-+u_0+u_+,u_-\tf+u_0\tf+u_+\tf)
&=\sigma_-\langle  u_-,u_-\tf \rangle_{H^1_0(\domL)}
+\sum_{m\in\setN} \dmL \beta_m\beta_m\tf
+\sigma_+\langle  u_+,u_+\tf \rangle_{H^1_0(\domL)}.
\end{align*}
\end{lemma}
\begin{proof}
The statement can be obtained by repeating the steps from \Cref{lem:blockoperator_semi} one-to-one. 
\end{proof}
We will now analyze the diagonal entries \(\dmL\) and to do this we introduce the following functions:
\begin{align*}
\tilde\mf_{\contrast,\ratio,\ev_{m}}(t)&:=
\frac{1+\frac{\contrast\sqrt{\ratio^2t^{2}+12}}{\sqrt{t^{2}+12}} \mz_{\ratio,\ev_{m}}(t)}%
{1+\frac{\sqrt{\ratio^2t^{2}+12}}{\sqrt{t^{2}+12}} \mz_{\ratio,\ev_{m}}(t)},&
\mz_{\ratio,\ev_{m}}(t)&:=\frac{\mj_{\ev_{m} L}(\mq(t))}{\mj_{\ev_{m} L}(\mq(\ratio t))},\\
\mq(t)&:=\left(\frac{1+\frac{1}{3}t^{2}-t\sqrt{1+\frac{1}{12}t^{2}}}{1+\frac{1}{3}t^{2}+t\sqrt{1+\frac{1}{12}t^{2}}}\right)^{\frac{1}{t}},&
\mj_{n}(q)&:=\frac{1-q^{n}}{1+q^{n}}.
\end{align*}
\begin{lemma}\label{lem:diagonal_L}
The diagonal entries $\dmL$ defined in \eqref{eq:dm_L} satisfy 
$\dmL=\sigma_-\tilde\mf_{\contrast,\ratio,\ev_m}(\ev_m\hxm)$.
\end{lemma}
\begin{proof}
We first compute
\begin{align*}
\nu_{m,-}^{\Nxm}&=\left(\frac{\mummI}{\mummII}\right)^{\frac{L}{\hxm}}
=\left(\frac{1+\frac{1}{3}\ev_{m}^2\hxm^{2}-\hxm\ev_{m}\sqrt{1+\frac{1}{12}\ev_{m}^{2}\hxm^{2}}}{1+\frac{1}{3}\ev_{m}^2\hxm^{2}+\hxm\ev_{m}\sqrt{1+\frac{1}{12}\ev_{m}^{2}\hxm^{2}}}\right)^{\frac{\ev_{m}L}{\ev_{m}\hxm}}
=\mq(\hxm\ev_{m})^{\ev_{m}L}.
\end{align*}
In the same way, we get \(\nu_{m,+}^{\Nxp}=\mq(r\hxm\ev_{m})^{\ev_{m}L}\).
Now we calculate
\begin{align*}
\dmL&=
\frac{\sigma_{-}\frac{1+\nu_{m,-}^{\Nxm}}{1-\nu_{m,-}^{\Nxm}}\sqrt{(\amm)^2-(\bmm)^2}+\sigma_{+}\frac{1+\nu_{m,+}^{\Nxp}}{1-\nu_{m,+}^{\Nxp}}\sqrt{(\amp)^2-(\bmp)^2}}%
{\phantom{\sigma_{-}}\frac{1+\nu_{m,-}^{\Nxm}}{1-\nu_{m,-}^{\Nxm}}\sqrt{(\amm)^2-(\bmm)^2}+\phantom{\sigma_{+}}\frac{1+\nu_{m,+}^{\Nxp}}{1-\nu_{m,+}^{\Nxp}}\sqrt{(\amp)^2-(\bmp)^2}}\\
&=\frac{\frac{\sigma_{-}\ev_{m}}{\mj_{\ev_{m} L}(\mq(\ev_{m}\hxm))}\sqrt{1+\frac{1}{12}\ev_{m}^{2}\hxm^{2}}+\frac{\sigma_{+}\ev_{m}}{\mj_{\ev_{m} L}(\mq(\ratio\ev_{m}\hxm))}\sqrt{1+\frac{1}{12}\ratio^{2}\ev_{m}^{2}\hxm^{2}}}%
{\phantom{\sigma_{-}}\frac{\ev_{m}}{\mj_{\ev_{m} L}(\mq(\ev_{m}\hxm))}\sqrt{1+\frac{1}{12}\ev_{m}^{2}\hxm^{2}}+\phantom{\sigma_{+}}\frac{\ev_{m}}{\mj_{\ev_{m} L}(\mq(\ratio\ev_{m}\hxm))}\sqrt{1+\frac{1}{12}\ratio^{2}\ev_{m}^{2}\hxm^{2}}}\\
&=\sigma_{-}\frac{\sqrt{1+\frac{1}{12}\ev_{m}^{2}\hxm^{2}}+\frac{\contrast\mj_{\ev_{m} L}(\mq(\ev_{m}\hxm))}{\mj_{\ev_{m} L}(\mq(\ratio\ev_{m}\hxm))}\sqrt{1+\frac{1}{12}\ratio^{2}\ev_{m}^{2}\hxm^{2}}}%
{\sqrt{1+\frac{1}{12}\ev_{m}^{2}\hxm^{2}}+\frac{\mj_{\ev_{m} L}(\mq(\ev_{m}\hxm))}{\mj_{\ev_{m} L}(\mq(\ratio\ev_{m}\hxm))}\sqrt{1+\frac{1}{12}\ratio^{2}\ev_{m}^{2}\hxm^{2}}}\\
&=\sigma_{-}\frac{1+\contrast\mz_{\ratio,\ev_{m}}(\ev_{m}\hxm)\frac{\sqrt{12+\ratio^{2}\ev_{m}^{2}\hxm^{2}}}{\sqrt{12+\ev_{m}^{2}\hxm^{2}}}}%
{1+\mz_{\ratio,\ev_{m}}(\ev_{m}\hxm)\frac{\sqrt{12+\ratio^{2}\ev_{m}^{2}\hxm^{2}}}{\sqrt{12+\ev_{m}^{2}\hxm^{2}}}}
=\sigma_{-}\tilde\mf_{\contrast,\ratio,\ev_m}(\ev_m\hxm).
\end{align*}
\end{proof}
In preperation of \Cref{lem:lim_f_zero_L} we formulate the following Lemma.
\begin{lemma}\label{lem:lim_z}
For each $c>0$ it holds that
$\lim_{\hxm\to0+} \sup_{\ev>0,\ev\hxm\geq c} \mz_{\ratio,\ev}(\ev\hxm)=1$.
\end{lemma}
\begin{proof}
Note that the auxiliary function $\mj_n$ in the definition of $\mz_{\ratio,\ev_m}$ was chosen this way to deal with a particular technicality in proof \Cref{lem:liminf_f_L}.
Here however, we exploit the more explicit representation
\begin{align*}
    \mz_{\ratio,\ev}&(\ev\hxm)=\frac{1-\left(\mq(\ev\hxm)^{\ev\hxm}\right)^{\frac{L}{\hxm}}}{1+\left(\mq(\ev\hxm)^{\ev\hxm}\right)^{\frac{L}{\hxm}}}\frac{1+\left(\mq(\ratio\ev\hxm)^{\ratio\ev\hxm}\right)^{\frac{L}{\ratio\hxm}}}{1-\left(\mq(\ratio\ev\hxm)^{\ratio\ev\hxm}\right)^{\frac{L}{\ratio\hxm}}}
\end{align*}
The claim follows now from
$\sup_{t\geq\frac{c}{\max(1,\ratio)}}\mq(t)^t<1$.
\end{proof}

The following lemma is the pendant to \Cref{lem:lim_f_zero}.
\begin{lemma}\label{lem:lim_f_zero_L}
If \eqref{eq:AssR} is satisfied, then $\lim_{\hxm\to0+} \inf_{m\in\setN} \big|\tilde\mf_{\contrast,\ratio,\ev_{m}}(\ev_m\hxm)\big|=0$.
\end{lemma}
\begin{proof}
We proceed as in the proof of \Cref{lem:lim_f_zero}:

Let $\hxm=\frac{1}{l+\epsilon} \roott$ with $l\in\setN_0$ and $\epsilon\in[0,1)$.
We choose $m=l$ and exploit that $\lim_{l\to+\infty}\frac{l}{l+\epsilon}\to1$ uniformly in $\epsilon\in[0,1)$.
Thus also $\ev_m\hxm=\frac{l}{l+\epsilon}\roott\xrightarrow{l\to+\infty}\roott$ uniformly in $\epsilon\in[0,1)$.

\Cref{lem:lim_z} yields that $\lim_{l\to+\infty}\mz_{\ratio,\ev_l}(\frac{l}{l+\epsilon}\roott)=1$, which in combination with the continuity of $\mf_{\contrast,\ratio}$ and $\mf_{\contrast,\ratio}(\roott)=0$ provides the claim.
\end{proof}

In preperation of \Cref{lem:liminf_f_L} we introduce the following \Cref{lem:func_properties,lem:cont_z}, where \Cref{lem:func_properties} is itself an auxiliary result for \Cref{lem:cont_z}.
\begin{lemma}\label{lem:func_properties}
Let \(I\) be a closed subintervall of \((0,1)\) and \(c>0\).
Then the following statements hold:
\begin{enumerate}
\item The family \(\{\mj_{n}:n\in[c,\infty)\}\) and the family of their derivatives are both uniformly bounded on \(I\).
\item The family \(\{\frac{1}{\mj_{n}}:n\in[c,\infty)\}\) and the family of their derivatives are both uniformly bounded on \(I\).
\item \(\mq\) is continuously differentiable on \([0,\infty)\).
\item \(\lim_{t\rightarrow 0^{+}}\mq(t)=e^{-2}\).
\end{enumerate}
\end{lemma}
\begin{proof}

Let $I=[a,b]$, $a,b\in(0,1)$.
Thence $q^n\in (0,b^c]\subset(0,1)$ for $q\in[a,b]$.
Thus $\mj_{n}(q)\leq 1$ and $1/\mj_{n}(q)\leq\sup_{s\in(0,b^c]} \frac{1+s}{1-s}=\frac{1+b^c}{1-b^c}$.
Furthermore, $\mj_{n}'(q)=\frac{2nq^{n-1}}{(1-q^n)^2}$ and $(1/\mj_{n})'(q)=\frac{-2nq^{n-1}}{(1+q^n)^2}$ from which we can deduce the first two claims.
The differentiability of $\mq$ follows in a straightforward fashion.
The last claim follows by applying the de l'Hospital rule to $\log\mq$.
\end{proof}
\begin{lemma}\label{lem:cont_z}
For \(\ratio,c>0\) the family \(\{\mz_{\ratio,\ev}:\ev\in[c,\infty)\}\) is equicontinuous from the right at \(0\) and \(\lim_{t\rightarrow 0+}\mz_{\ratio,\ev}(t)=1\).
\end{lemma}
\begin{proof}
Due to the second half of \Cref{lem:func_properties} there exist $\delta>0$ and $q_1,q_2\in(0,1),q_1<q_2$ such that $\mq(t),\mq(rt)\in(q_1,q_2)$ for all $t\in[0,\delta)$.
Then \Cref{lem:func_properties} and the product rule yield that $t\mapsto\mz_{\ratio,\ev}(t)$ is uniformly equicontinuous on $t\in[0,\delta)$, $\ev\geq c$.
At last we obtain \(\lim_{t\rightarrow 0+}\mz_{\ratio,\ev}(t)=1\) from \(\lim_{t\rightarrow 0+}\mq(t)=e^{-2}=\lim_{t\rightarrow 0+}\mq(\ratio t)\).
\end{proof}
\begin{lemma}\label{lem:liminf_f_L}
If $\epsilon\in(0,1)$ and one of the following two conditions
\begin{align}
\label{eq:AssRL_vv}
\underbrace{|\contrast|(1+\epsilon)<1 \text{ and } \ratio|\contrast|(1+\epsilon)<1
}_{(\ref{eq:AssRL_vv}a)}
\quad\text{ or }\quad
\underbrace{|\contrast|(1-\epsilon)>1 \text{ and } \ratio|\contrast|(1-\epsilon)>1}_{(\ref{eq:AssRL_vv}b)}
\end{align}
is satisfied, then there exists \(\delta>0\) such that
\[\inf_{\substack{\hxm\in(0,\delta)\\\lambda\geq1}} |\tilde\mf_{\contrast,\ratio,\ev}(\lambda\hxm)|\geq
\min_{p\in\{\pm1\}}\min\left\{\frac{|1+(1+p\epsilon)\contrast|}{2+\epsilon},\frac{|1+(1+p\epsilon)\contrast\ratio|}{1+(1+\epsilon)\ratio}\right\}.
\]
\end{lemma}
\begin{proof}
Due to \Cref{lem:cont_z} we can find \(\tau>0\) such that
\(\mz_{\ratio,\ev}(\ev\hxm)\in[1-\epsilon,1+\epsilon]\) for all \(\ev\hxm\in[0,\tau], \hxm>0, \ev\geq1\).
On the other hand, \Cref{lem:lim_z} yields the existence of $\delta>0$ such that
\(\mz_{\ratio,\ev}(\ev\hxm)\in[1-\epsilon,1+\epsilon]\) for all \(\ev\hxm\geq\tau, \hxm\in(0,\delta), \ev\geq1\).
Thus \(\mz_{\ratio,\ev}(t)\in[1-\epsilon,1+\epsilon]\) for all \(t=\ev\hxm, \ev\geq1, \hxm\in(0,\delta)\).
Thence
\[\inf_{\substack{\hxm\in(0,\delta)\\\lambda\geq1}} |\tilde\mf_{\contrast,\ratio,\ev}(\lambda\hxm)|
\geq \inf_{t>0,p\in\{\pm1\}}
\frac{\left|1+\contrast(1+p\epsilon)\sqrt{\frac{12+\ratio^{2}t^{2}}{12+t^{2}}}\right|}%
{1+(1+\epsilon)\sqrt{\frac{12+t^{2}\hxm^{2}}{12+t^{2}}}},
\]
from which the claim follows.
\end{proof}
Having analyzed \(\tilde\mf_{\contrast,\ratio,\ev_{m}}\) we can now proof the following theorem about the stability of \(\opL_{\hxpm}\).
\begin{theorem}\label{thm:semi_L}
If for some \(\epsilon\in(0,1)\) \eqref{eq:AssRL_vv} is satisfied, then $\opL_{\hxpm}^{-1}$ exists and satisfies
\begin{align*}
\|\opL_{\hxpm}^{-1}\|_{\Lspace(\VsemiL)}\leq
\frac{1}{|\sigma_{-}|\min\left\{1,|\kappa|,\frac{|1+(1+\epsilon)\contrast|}{2+\epsilon},\frac{|1+(1-\epsilon)\contrast|}{2+\epsilon},\frac{|1+(1+\epsilon)\contrast\ratio|}{1+(1+\epsilon)\ratio},\frac{|1+(1-\epsilon)\contrast\ratio|}{1+(1+\epsilon)\ratio}\right\}}
\end{align*}
for all $\hxm\in(0,\delta)$ with $\delta>0$ as in \Cref{lem:liminf_f_L}.
Contrary, if \eqref{eq:AssR} is satisfied then
\begin{align*}
\lim_{\hxm\rightarrow0} \|\opL_{\hxpm}^{-1}\|_{\Lspace(\VsemiL)} =+\infty
\end{align*}
(where we define $\|\opL_{\hxpm}^{-1}\|_{\Lspace(\VsemiL)}:=+\infty$, if $\opL_{\hxpm}^{-1}$ does not exist).
\end{theorem}
\begin{proof}
As in the previous sections the theorem follows directly from the properties of \(\tilde\mf_{\contrast,\ratio,\ev_{m}}\) that where shown in the Lemmas \ref{lem:lim_f_zero_L} and \ref{lem:liminf_f_L}.
\end{proof}

\subsubsection{Full discretization}

As for the unbounded domain, the only difference is that we now consider \(\evdism\) that also depend on \(\hxm\) and where in variables that depend on \(\ev_{m}\) we replace it by \(\evdism\) and indicate this by adding a hat.
In addition, we define \(\hat{\nu}_{m,\pm}:=\frac{\mumpmIhat}{\mumpmIIhat}\).
The following three lemmas can then be derived correspondingly to \Cref{lem:decomposition_semi_L,lem:blockoperator_semi_L,lem:diagonal_L}.
\begin{lemma}\label{lem:decomposition_full_L}
The space $\VfullL$ admits an orthogonal decomposition
\begin{align*}
\VfullL
=\Big(\VxL^-\otimes\Vy\Big)
\oplus^\bot \VxzhatL
\oplus^\bot \Big(\VxL^+\otimes\Vy\Big),
\end{align*}
where $\VxL^-:=\{v\in\VxL\colon v|_{\setR_+}=0\}$ and
$\VxL^+-:=\{v\in\VxL\colon v|_{\setR_-}=0\}$.
The subspace $\VxzhatL$ is spanned by the orthonormal basis
$\big(\vmhatL(x)\otimes \bfydis_m(y)\big)_{m=1,\dots,\Ny-1}$, where
\begin{align*}
\vmhatL(x)&:=\frac{1}{\sqrt{\frac{1+\hat{\nu}_{m,-}^{\Nxm}}{1-\hat{\nu}_{m,-}^{\Nxm}}\sqrt{(\ammhat)^2-(\bmmhat)^2}+\frac{1+\hat{\nu}_{m,+}^{\Nxp}}{1-\hat{\nu}_{m,+}^{\Nxp}}\sqrt{(\amphat)^2-(\bmphat)^2}}}
\Bigg(\hfx_0(x)\\
&\quad+\sum_{n=-\Nxm}^{-1}\frac{\mummIhat^{-n}-\hat{\nu}_{m,-}^{\Nxm}\mummIIhat^{- n}}{1-\hat{\nu}_{m,-}^{\Nxm}}\hfx_n(x)
+\sum_{n=1}^{\Nxp}\frac{\mumpIhat^{n}-\hat{\nu}_{m,+}^{\Nxp}\mumpIIhat^{n}}{1-\hat{\nu}_{m,+}^{\Nxp}}\hfx_n(x)\Bigg).
\end{align*}
\end{lemma}
\begin{lemma}\label{lem:blockoperator_full_L}
The operator $\opL_{\hxpm,\hy}$ is block diagonal with respect to the orthogonal decomposition given in \Cref{lem:decomposition_full_L}.
The blocks corresponding to $\VxL^-\otimes\Vy$ and $\VxL^+\otimes\Vy$ equal the identity times $\sigma_-$ and $\sigma_+$ respectively.
The block corresponding to $\VxzhatL$ is diagonal with respect to the basis given in \Cref{lem:decomposition_full_L} and the diagonal entries are given by
\begin{align}\label{eq:dm_L_full}
\dmhatL:=
\frac{\sigma_{-}\frac{1+\hat\nu_{m,-}^{\Nxm}}{1-\hat\nu_{m,-}^{\Nxm}}\sqrt{(\amm)^2-(\bmm)^2}+\sigma_{+}\frac{1+\hat\nu_{m,+}^{\Nxp}}{1-\hat\nu_{m,+}^{\Nxp}}\sqrt{(\amp)^2-(\bmp)^2}}%
{\phantom{\sigma_{-}}\frac{1+\hat\nu_{m,-}^{\Nxm}}{1-\hat\nu_{m,-}^{\Nxm}}\sqrt{(\amm)^2-(\bmm)^2}+\phantom{\sigma_{+}}\frac{1+\hat\nu_{m,+}^{\Nxp}}{1-\hat\nu_{m,+}^{\Nxp}}\sqrt{(\amp)^2-(\bmp)^2}},
\quad m=1,\dots,\Ny-1.
\end{align}
Non-verbally: For each
$u_-,u_-\tf\in\VxL^-\otimes\Vy$,
$u_0,u_0\tf\in\VxzhatL$,
$u_+,u_+\tf\in\VxL^+\otimes\Vy$
and $u_0=\sum_{m=1}^{\Ny-1} \beta_m \vmhatL\otimes\bfydis_m$,
$u_0\tf=\sum_{m=1}^{\Ny-1} \beta_m\tf \vmhatL\otimes\bfydis_m$,
$(\beta_m)_{m=1}^{\Ny-1},(\beta_m\tf)_{m=1}^{\Ny}\in\setR^{\Ny-1}$
it holds that
\begin{align*}
\ses_{\domL}(u_-+u_0+u_+,u_-\tf+u_0\tf+u_+\tf)
&=\sigma_-\langle  u_-,u_-\tf \rangle_{H^1_0(\domL)}
+\sum_{m=1}^{\Ny-1} \dmhatL \beta_m\beta_m\tf
+\sigma_+\langle  u_+,u_+\tf \rangle_{H^1_0(\domL)}.
\end{align*}
\end{lemma}
\begin{lemma}\label{lem:dmhatL}
The diagonal entries $\dmhatL$ satisfy
$\dmhatL=\sigma_-\tilde\mf_{\contrast,\ratio,\evdism}\big(\mh_{\ratioy}(\ratioy m\hxm)\big)$
for $m=1,\dots,\Ny-1$.
\end{lemma}
Now we investigate \(\tilde\mf_{\contrast,\ratio,\evdism}\big(\mh_{\ratioy}(\ratioy m\hxm)\big)\).
Since we have already shown, that \(\tilde\mf_{\contrast,\ratio,\ev}\) is equicontinuous at zero in \(\ev\geq1\) and because \(\evdism\geq\ev_m\geq1\),
the only difference to the previous analysis concerning the unbounded domain is that
we have to deal with the additional composition with the continuous function \(\mh_{\ratioy}\).
Before we formulate the next \Cref{lem:lim_fh_zero}, let us recall that $\hxm=\frac{\hy}{\ratioy}=\frac{\pi}{\ratioy\Ny}$.
\begin{lemma}\label{lem:lim_fh_zero_L}
If \eqref{eq:AssRy} is satisfied, then $\lim_{\Ny\to+\infty} \inf_{m\in\{1,\dots,\Ny-1\}}
\big|\tilde\mf_{\contrast,\ratio,\evdism}\big(\mh_{\ratioy}(\frac{m\pi}{\ratioy\Ny})\big)\big|=0$.
\end{lemma}
\begin{proof}
The proof follows along the lines of the proof of \Cref{lem:lim_fh_zero}, where in addition we apply \Cref{lem:lim_z} to cope with the replacement of $\mf_{\contrast,\ratio}$ by $\tilde\mf_{\contrast,\ratio,\evdism}$.
\end{proof}
\begin{lemma}\label{lem:liminf_f_full_L}
If for \(\epsilon \in (0,1)\) one of the following two conditions
\begin{subequations}\label{eq:AssRyL_vv}
\begin{align}
|\contrast|(1+\epsilon)<1 \,&\text{ and }\, \ratio^2\contrast^2(1+\epsilon)^2<1+\ratioy(1-\contrast^2(1+\epsilon)^2) \\
&\text{ or }\nonumber\\
|\contrast|(1-\epsilon)>1 \,&\text{ and }\, \ratio^2\contrast^2(1-\epsilon)^2>1+\ratioy^2(1-\contrast^2(1-\epsilon)^2) 
\end{align}
\end{subequations}
is satisfied, then there exists \(\delta>0\) such that
\begin{align*}
\inf_{\substack{\hxm\in(0,\delta)\\m\in\{1,\dots,\Ny-1\}}}
\big|\tilde\mf_{\contrast,\ratio,\evdism}&\big(\mh_{\ratioy}(\ratioy m\hxm)\big)\big|\\
&\geq
\min_{p\in\{\pm1\}}
\min\Bigg\{\frac{|1+(1+p\epsilon)\contrast|}{2+\epsilon},
\frac{\big|\sqrt{1+\ratioy^2}+(1+p\epsilon)\contrast\sqrt{\ratio^2+\ratioy^2}\big|}{\sqrt{1+\ratioy^2}+(1+\epsilon)\sqrt{\ratio^2+\ratioy^2}}\Bigg\}>0.
\end{align*}
\end{lemma}
\begin{proof}
It suffices to combine the techniques used for \Cref{lem:liminf_f_full} and \Cref{lem:liminf_f_L}. As in the proof of \Cref{lem:liminf_f_L} we choose \(\delta>0\) such that
\(\mz_{\ratio,\evdism}(\evdism\hxm)\in [1-\epsilon,1+\epsilon]\) for all $\evdism\geq1, \hxm\in(0,\delta)$.
Thence
\begin{align*}
\inf_{\substack{\hxm\in(0,\delta)\\m\in\{1,\dots,\Ny-1\}}}
\big|\tilde\mf_{\contrast,\ratio,\evdism}&\big(\mh_{\ratioy}(\ratioy m\hxm)\big)\big|
\geq
\inf_{t\in[0,\sqrt{12}/\ratioy],p\in\{\pm1\}}
\frac{\left|1+\contrast(1+p\epsilon)\sqrt{\frac{12+\ratio^{2}t^{2}}{12+t^{2}}}\right|}%
{1+(1+\epsilon)\sqrt{\frac{12+t^{2}\hxm^{2}}{12+t^{2}}}},
\end{align*}
from which the claim follows.
\end{proof}
Now we are in the position to conclude our analysis of the full discretization of \eqref{al:org_prob_dir} in the following theorem.
\begin{theorem}\label{thm:full_L}
If \eqref{eq:AssRyL_vv} is satisfied for some \(\epsilon\in(0,1)\), then $\opL_{\hxpm,\hy}^{-1}$ exists and satisfies
\begin{align*}
\|\opL_{\hxpm,\hy}^{-1}\|_{\Lspace(\VfullL)}
\leq
\frac{1}{|\sigma_{-}|\min_{p\in\{\pm1\}}\min\left\{1,|\kappa|,\frac{|1+(1+p\epsilon)\contrast|}{2+\epsilon},\frac{|\sqrt{1+\ratioy^2}+(1+p\epsilon)\contrast\sqrt{\ratio^2+\ratioy^2}|}{\sqrt{1+\ratioy^2}+(1+\epsilon)\sqrt{\ratio^2+\ratioy^2}}\right\}}
\end{align*}
for all $\hxm\in(0,\delta)$ with $\delta>0$ as in \Cref{lem:liminf_f_full}. Contrary, if \eqref{eq:AssRy} is satisfied, then
\begin{align*}
\lim_{\hxm\rightarrow0} \|\opL_{\hxpm,\hy}^{-1}\|_{\Lspace(\VfullL)}=+\infty
\end{align*}
(where we define $\|\opL_{\hxpm,\hy}^{-1}\|_{\Lspace(\VfullL)}:=+\infty$, if $\opL_{\hxpm,\hy}^{-1}$ does not exist).
\end{theorem}
\begin{proof}
As for \Cref{thm:semi,thm:full,thm:semi_L} the claims follow directly from the respective \Cref{lem:blockoperator_full_L,lem:dmhatL,lem:lim_fh_zero_L,lem:liminf_f_full_L}.
\end{proof}

\section{Computational examples}\label{sec:num}
We will now confirm our theoretical results by testing them in explicit computational examples. The code to reproduce all of them is provided in \cite{BHHLLP:25}.
We consider the following problem posed on a bounded domain:
\begin{align*}
\begin{aligned}
	\text{Find } u \in H^1_0(\domL) \text{ such that }
	-\div(\sigma\nabla u)=f \text{ in }\domL=(-L,L)\times(0,\pi),\\
 \text{with } f(x,y)=-\sigma\left[\frac{-2y(y-\pi)(y-2\pi)}{L^2}+6\left(1-\frac{x^2}{L^2}\right)(y-\pi)\right],
\end{aligned}
\end{align*}
which has the solution
\[u(x,y)=\left(1-\frac{x^{2}}{L^{2}}\right)y(y-\pi)(y-2\pi).\] 

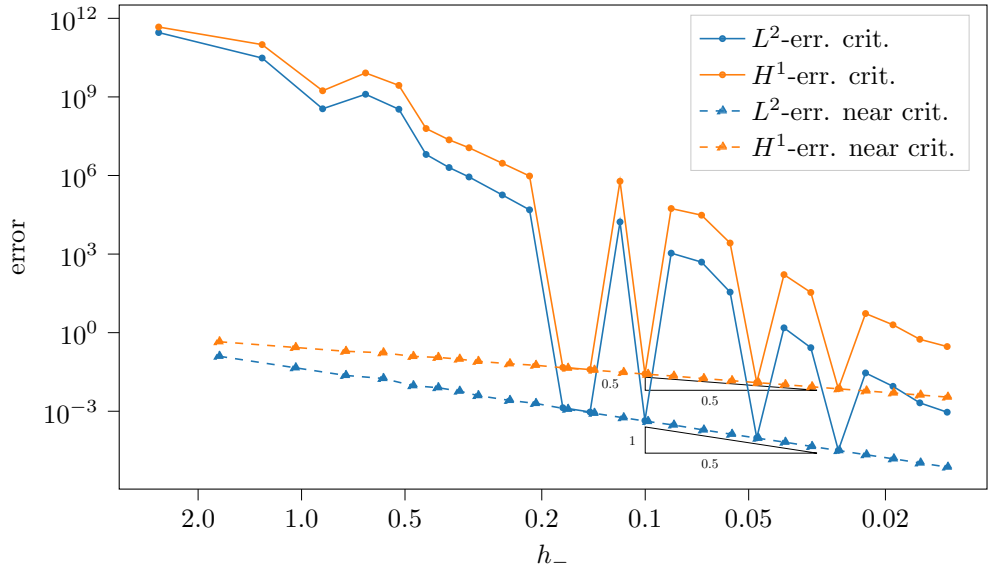
\begin{figure}
\centering
\begin{tikzpicture}[scale=1]

\definecolor{darkgray176}{RGB}{176,176,176}
\definecolor{darkorange25512714}{RGB}{255,127,14}
\definecolor{lightgray204}{RGB}{204,204,204}
\definecolor{steelblue31119180}{RGB}{31,119,180}

\begin{axis}[width=13cm,height=8cm,
legend cell align={left},
legend style={fill opacity=0.8, draw opacity=1, text opacity=1, draw=lightgray204},
log basis x={10},
log basis y={10},
tick align=outside,
tick pos=left,
x grid style={darkgray176},
xlabel={\(\displaystyle h_{-}\)},
xmin=0.294737517371598, xmax=98.7550313071886,
xmode=log,
xtick style={color=black},
xtick={0.5,1,2,5,10,20,50,100},
xticklabels={2.0,1.0,0.5,0.2,0.1,0.05,0.02,0.01},
y grid style={darkgray176},
ylabel={error},
ymin=1.06591727242928e-06, ymax=3221648228451.59,
ymode=log,
ytick style={color=black},
ytick={1e-09,1e-06,0.001,1,1000,1000000,1000000000,1000000000000,1e+15,1e+18},
yticklabels={
  \(\displaystyle {10^{-9}}\),
  \(\displaystyle {10^{-6}}\),
  \(\displaystyle {10^{-3}}\),
  \(\displaystyle {10^{0}}\),
  \(\displaystyle {10^{3}}\),
  \(\displaystyle {10^{6}}\),
  \(\displaystyle {10^{9}}\),
  \(\displaystyle {10^{12}}\),
  \(\displaystyle {10^{15}}\),
  \(\displaystyle {10^{18}}\)
}
]
\path [draw=black, very thin]
(axis cs:10,0.000251188643150958)
--(axis cs:31.6227766016838,2.51188643150958e-05)
--(axis cs:10,2.51188643150958e-05)
--cycle;
\path [draw=black, very thin]
(axis cs:10,0.0199526231496888)
--(axis cs:31.6227766016838,0.00630957344480193)
--(axis cs:10,0.00630957344480193)
--cycle;
\addplot [semithick, steelblue31119180, mark=*, mark size=1, mark options={solid}]
table {%
0.383896167102679 286245228391.311
0.767792334205357 30501296361.8988
1.15168850130804 351882224.993559
1.53558466841071 1262890933.61168
1.91948083551339 338210013.39803
2.30337700261607 6353453.28178043
2.68727316971875 2013259.57623923
3.07116933682143 879421.61492825
3.83896167102679 180694.558475065
4.60675400523214 48905.7770356016
5.75844250654018 0.00136095478338602
6.91013100784821 0.000934389801132374
8.44571567625893 16864.1397372964
9.98130034466964 0.000440103190728086
11.900781180183 1086.54338215455
14.5880543499018 493.467442474987
17.6592236867232 35.3046465190135
21.1142891906473 9.63535897514613e-05
25.3371470287768 1.52869214445269
30.3277972011116 0.267070860852725
36.4701358747545 3.20482506164159e-05
43.7641630497054 0.0289998534800838
52.593774893067 0.0089377528449195
62.9589714048393 0.00208964621109557
75.6275449192277 0.000919600805172264
};
\addlegendentry{$L^2$-err. crit.}
\addplot [semithick, darkorange25512714, mark=*, mark size=1, mark options={solid}]
table {%
0.383896167102679 465652564469.209
0.767792334205357 99236636714.3851
1.15168850130804 1717284805.54444
1.53558466841071 8217686704.30918
1.91948083551339 2750934241.82324
2.30337700261607 62013298.8164537
2.68727316971875 22925644.6629668
3.07116933682143 11444870.0525545
3.83896167102679 2939468.54668713
4.60675400523214 954696.109346587
5.75844250654018 0.0470901777375446
6.91013100784821 0.0390190075914719
8.44571567625893 603546.333729891
9.98130034466964 0.0267788896676739
11.900781180183 54793.9372669954
14.5880543499018 30504.6358163383
17.6592236867232 2641.88216925669
21.1142891906473 0.0125299438181903
25.3371470287768 164.129899956196
30.3277972011116 34.3223742771102
36.4701358747545 0.00722630422251363
43.7641630497054 5.37804716764832
52.593774893067 1.99192631518535
62.9589714048393 0.557502845618136
75.6275449192277 0.29471734028125
};
\addlegendentry{$H^1$-err. crit.}
\addplot [semithick, steelblue31119180, dashed,mark=triangle*, mark size=2]
table {%
0.575844250654018 0.125633414048247
0.959740417756696 0.0455884233023563
1.34363658485937 0.0233102719097896
1.72753275196205 0.0178507793027766
2.11142891906473 0.00945251734458609
2.49532508616741 0.00794227438102325
2.87922125327009 0.00583651886078304
3.26311742037277 0.0039597588446747
4.03090975457812 0.00259529113613352
4.79870208878348 0.00198693201159898
5.95039059009152 0.00119113528404166
7.10207909139955 0.000836173022644225
8.63766375981027 0.000565309278651527
10.173248428221 0.000407536137872354
12.0927292637344 0.000297797691340472
14.7800024334531 0.000193083673927701
17.8511717702746 0.000132362114093948
21.3062372741987 9.29148897736777e-05
25.5290951123281 6.57016240143484e-05
30.5197452846629 4.52834028929877e-05
36.6620839583058 3.13810276722674e-05
43.9561111332567 2.18304826970203e-05
52.7857229766183 1.52485784558241e-05
63.1509194883906 1.05765246326929e-05
75.819493002779 7.37461952198657e-06
};
\addlegendentry{$L^2$-err. near crit.}
\addplot [semithick, darkorange25512714, dashed,mark=triangle*, mark size=2]
table {%
0.575844250654018 0.450265171134009
0.959740417756696 0.272080263581882
1.34363658485937 0.19471945277066
1.72753275196205 0.170459499474192
2.11142891906473 0.124060616820629
2.49532508616741 0.113736842591611
2.87922125327009 0.0975063959004666
3.26311742037277 0.0803124316530251
4.03090975457812 0.0650224427130638
4.79870208878348 0.0568975682509929
5.95039059009152 0.0440527691513009
7.10207909139955 0.0369102008025899
8.63766375981027 0.030349092449146
10.173248428221 0.0257684460664507
12.0927292637344 0.022028066473948
14.7800024334531 0.01773704884193
17.8511717702746 0.0146855900323015
21.3062372741987 0.0123041846243224
25.5290951123281 0.0103467235099344
30.5197452846629 0.00858974888123319
36.6620839583058 0.00715063717492749
43.9561111332567 0.00596407344541846
52.7857229766183 0.00498457483222307
63.1509194883906 0.0041512891705034
75.819493002779 0.00346643262006986
};
\addlegendentry{$H^1$-err. near crit.}
\draw (axis cs:14.1253754462275,7.07945784384137e-06) node[
  scale=0.5,
  anchor=base west,
  text=black,
  rotate=0.0
]{0.5};
\draw (axis cs:8.7096358995608,5.01187233627273e-05) node[
  scale=0.5,
  anchor=base west,
  text=black,
  rotate=0.0
]{1};
\draw (axis cs:14.1253754462275,0.00177827941003892) node[
  scale=0.5,
  anchor=base west,
  text=black,
  rotate=0.0
]{0.5};
\draw (axis cs:7.2443596007499,0.00794328234724281) node[
  scale=0.5,
  anchor=base west,
  text=black,
  rotate=0.0
]{0.5};
\end{axis}

\end{tikzpicture}
\caption{Relative 
errors for critical values of \(\hxm\) (solid lines) and for nearly critical values of $\hxm$ (dashed lines).}
\label{fig:crit}
\end{figure}

\begin{figure}
\centering
\includegraphics[scale=0.2]{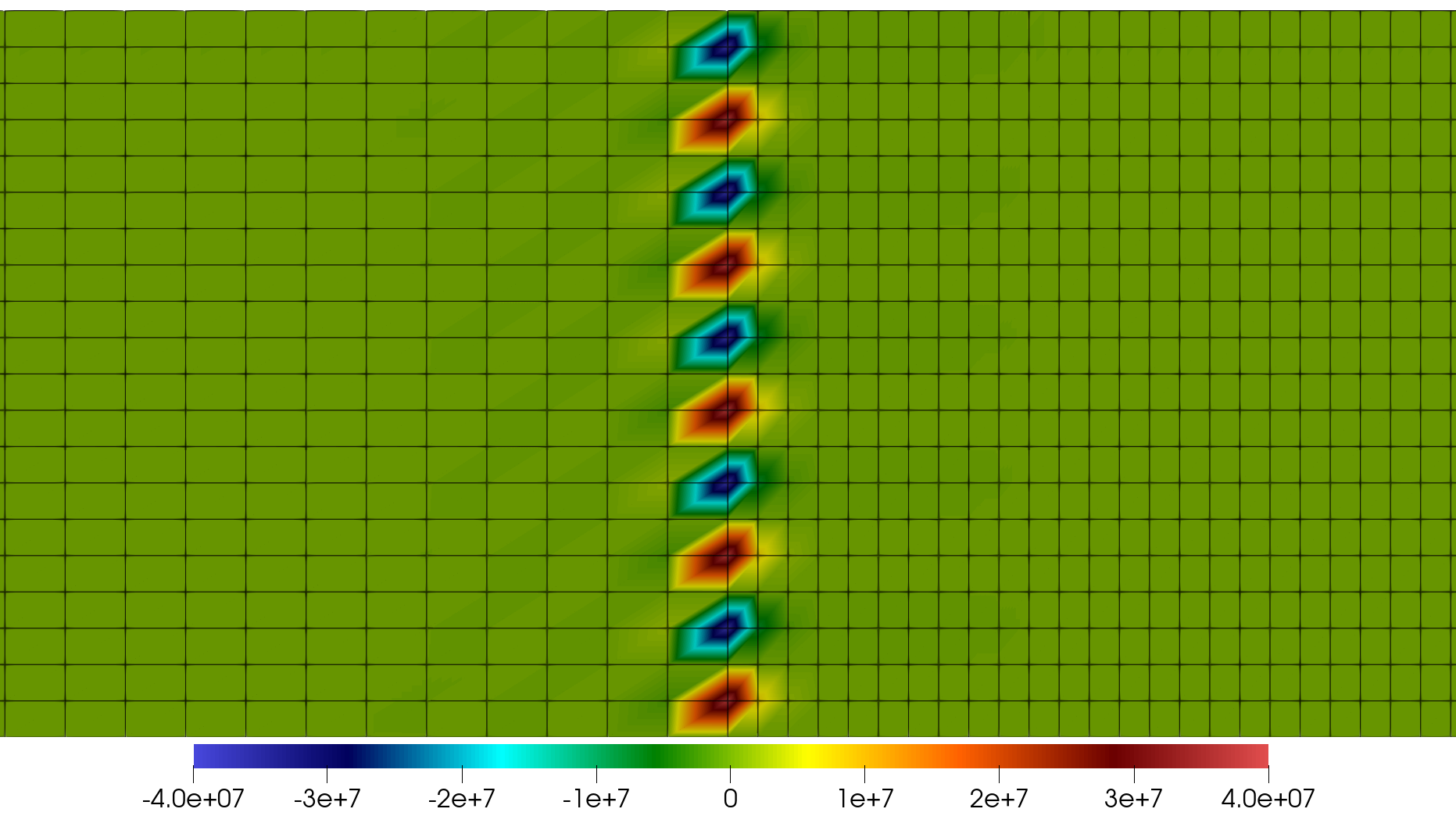}
\caption{
Computed solution for \(\hxm\approx0.26048\) in a neighborhood of the interface $\{0\}\times(0,\pi)$.}
\label{fig:crit_sol}
\end{figure}

First we examine the case of unstable discretizations.
To this end we consider parameters as follows:
\begin{align*}
	\sigma_{-}=-1,\quad
	\sigma_{+}=1.2,\quad
	\ratio=0.5,\quad
    \ratioy=\frac{2}{\sqrt{11}}.
\end{align*}
Hence the contrast at hand \(\kappa=-1.2\) is rather moderate.
Nevertheless, we are going to exhibit instabilities with this particular choice of parameters, whereas convenient examples of instabilities often require a much more critical contrast $|\contrast-\contrast_\mathrm{crit}|\approx 10^{-3}$ \cite[Fig.~1]{AbdulleHuberLemaire17}, \cite[Fig.~3]{BurmanPreussErn24}.
We choose a sufficiently large $L\approx26$ such that our problem is adequately close to an unbounded domain, and hence we can expect the critical values 
\begin{align}\label{eq:hcrit-example}
\hxm
:=\frac{1}{m}\frac{1}{\ratioy}\arccos\left(1+\frac{6\ratioy^{2}(1-\contrast^{2})}{(1-\contrast^{2}\ratio^{2})-2\ratioy^{2}(1-\contrast^{2})}\right)=\frac{\sqrt{11}}{2}\frac{1}{m}\arccos\left(1-1)\right)=\frac{\sqrt{11}\pi}{4m}\quad m\in\setN.
\end{align}
given in \Cref{thm:full} to yield also sensible values for our example.
In particular, we choose \(L:=10\frac{\sqrt{11}\pi}{4}\) such that
\begin{align*}
	\Nxm&=\frac{L}{\hxm}=10\frac{\sqrt{11}\pi}{4}\frac{4m}{\sqrt{11}\pi}=10m,\quad
	\Nxp=\frac{L}{\hxm}=10\frac{\sqrt{11}\pi}{4}\frac{8m}{\sqrt{11}\pi}=20m,\\
    \Ny&=\frac{\pi}{\ratioy\hxm}=\pi\frac{4m}{\sqrt{11}\pi}\frac{\sqrt{11}}{2}=2m
\end{align*}
are natural numbers of each \(m\in\setN\).
Even though we cannot expect our discretizations to have a non-trivial kernel at $\hxm$ we observe in \Cref{fig:crit} (solid lines) exorbitant errors.
Nevertheless, we recognize a decrease in the error which can be explained as follows:
The drastic error is triggered by the basis function $\vmhatL\otimes\bfydis_m$
for which the respective coefficient $\spl f,\vmhatL\otimes\bfydis_m\spr_{L^2}$ of $f$ decreases w.r.t.\ $m\in\setN$.
In \Cref{fig:crit_sol} we see the numerical solution for a critical value of \(\hxm\) and we observe 
the oscillating behaviour of $\vmhatL\otimes\bfydis_m$
that is corrupting the solution as predicted.

To further explore the possible behaviours of different discretizations we keep all parameters apart from $\hxm$ unchanged and choose now $\hxm=\frac{\sqrt{11}\pi}{4(m+\frac{1}{2})}$ to maximize the distance of $\hxm$ to the critical values \eqref{eq:hcrit-example} (while keeping $\Nxm,\Nxp,\Ny\in\setN$).
In contrast to our previous results, we observe in \Cref{fig:crit} (dashed lines) a distinct convergence of errors with convenient rates.

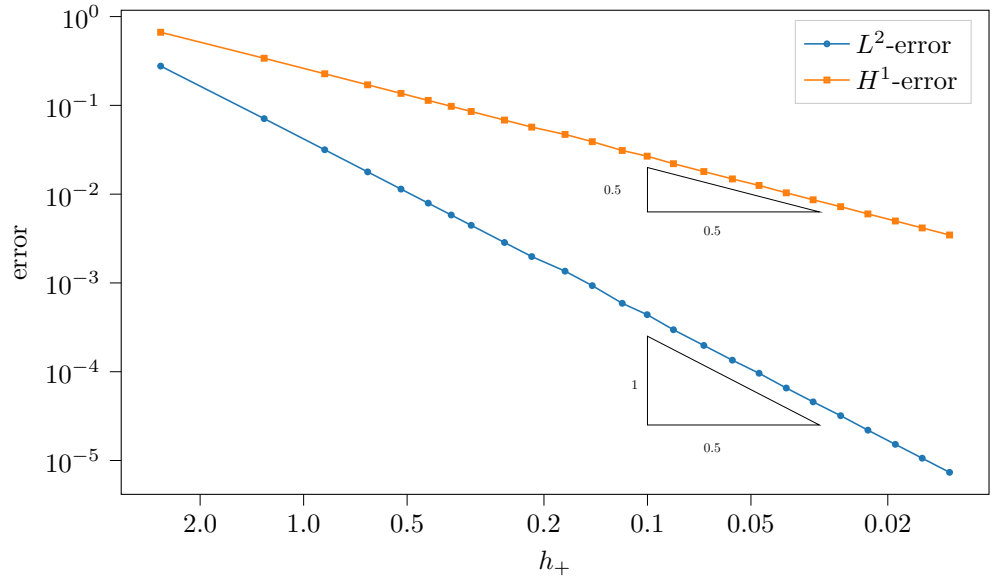
\begin{figure}
\centering
\begin{tikzpicture}

\definecolor{darkgray176}{RGB}{176,176,176}
\definecolor{darkorange25512714}{RGB}{255,127,14}
\definecolor{lightgray204}{RGB}{204,204,204}
\definecolor{steelblue31119180}{RGB}{31,119,180}

\begin{axis}[width=13cm,height=8cm,
legend cell align={left},
legend style={fill opacity=0.8, draw opacity=1, text opacity=1, draw=lightgray204},
log basis x={10},
log basis y={10},
tick align=outside,
tick pos=left,
x grid style={darkgray176},
xlabel={\(\displaystyle h_{+}\)},
xmin=0.294774875591359, xmax=98.4925345609346,
xmode=log,
xtick style={color=black},
xtick={0.5,1,2,5,10,20,50,100},
xticklabels={2.0,1.0,0.5,0.2,0.1,0.05,0.02,0.01},
y grid style={darkgray176},
ylabel={error},
ymin=4.15727086346968e-06, ymax=1.17826618168756,
ymode=log,
ytick style={color=black},
ytick={1e-07,1e-06,1e-05,0.0001,0.001,0.01,0.1,1,10,100},
yticklabels={
  \(\displaystyle {10^{-7}}\),
  \(\displaystyle {10^{-6}}\),
  \(\displaystyle {10^{-5}}\),
  \(\displaystyle {10^{-4}}\),
  \(\displaystyle {10^{-3}}\),
  \(\displaystyle {10^{-2}}\),
  \(\displaystyle {10^{-1}}\),
  \(\displaystyle {10^{0}}\),
  \(\displaystyle {10^{1}}\),
  \(\displaystyle {10^{2}}\)
}
]
\path [draw=black, very thin]
(axis cs:10,0.000251188643150958)
--(axis cs:31.6227766016838,2.51188643150958e-05)
--(axis cs:10,2.51188643150958e-05)
--cycle;
\path [draw=black, very thin]
(axis cs:10,0.0199526231496888)
--(axis cs:31.6227766016838,0.00630957344480193)
--(axis cs:10,0.00630957344480193)
--cycle;
\addplot [semithick, steelblue31119180, mark=*, mark size=1, mark options={solid}]
table {%
0.383896167102679 0.277064085892543
0.767792334205357 0.0708605177807353
1.15168850130804 0.0316205343946329
1.53558466841071 0.0178112390435849
1.91948083551339 0.0114064608363246
2.30337700261607 0.00792388507537095
2.68727316971875 0.00582283720249817
3.07116933682143 0.00445870874912759
3.83896167102679 0.00285402378244594
4.60675400523214 0.00198213055632657
5.75844250654018 0.00135774986785317
6.91013100784821 0.000932163775400765
8.44571567625893 0.000589805241943648
9.98130034466964 0.000439036068923166
11.900781180183 0.00029705922265868
14.5880543499018 0.000197698536675753
17.6592236867232 0.000134914196366848
21.1142891906473 9.61150875180212e-05
25.3371470287768 6.55372763602433e-05
30.3277972011116 4.57428181880531e-05
36.4701358747545 3.1968307520296e-05
43.7641630497054 2.19668814873127e-05
52.593774893067 1.52102818586473e-05
62.9589714048393 1.06142887007881e-05
75.6275449192277 7.35607723089434e-06
};
\addlegendentry{$L^2$-error}
\addplot [semithick, darkorange25512714, mark=square*, mark size=1, mark options={solid}]
table {%
0.383896167102679 0.665894540363034
0.767792334205357 0.339341555077619
1.15168850130804 0.227007517070679
1.53558466841071 0.170459813129422
1.91948083551339 0.136443379672831
2.30337700261607 0.113736990144196
2.68727316971875 0.097506505923429
3.07116933682143 0.0853282186843658
3.83896167102679 0.0682720060449979
4.60675400523214 0.056897607076178
5.75844250654018 0.0470901776997473
6.91013100784821 0.0390190075769985
8.44571567625893 0.031038775362715
9.98130034466964 0.0267788896655613
11.900781180183 0.0220280724863548
14.5880543499018 0.0179704231810168
17.6592236867232 0.0148452123276459
21.1142891906473 0.0125299438181463
25.3371470287768 0.0103467248522609
30.3277972011116 0.0086441139542196
36.4701358747545 0.0072263042225109
43.7641630497054 0.00599023156369713
52.593774893067 0.00498457514553016
62.9589714048393 0.00416394551686195
75.6275449192277 0.00346643277184047
};
\addlegendentry{$H^1$-error}
\draw (axis cs:14.1253754462275,1.25892541179417e-05) node[
  scale=0.5,
  anchor=base west,
  text=black,
  rotate=0.0
]{0.5};
\draw (axis cs:8.7096358995608,6.30957344480193e-05) node[
  scale=0.5,
  anchor=base west,
  text=black,
  rotate=0.0
]{1};
\draw (axis cs:14.1253754462275,0.00346736850452532) node[
  scale=0.5,
  anchor=base west,
  text=black,
  rotate=0.0
]{0.5};
\draw (axis cs:7.2443596007499,0.01) node[
  scale=0.5,
  anchor=base west,
  text=black,
  rotate=0.0
]{0.5};
\end{axis}

\end{tikzpicture}
\caption{Relative errors for a mesh satisfying \eqref{eq:AssRyL_vv}.}
\label{fig:uncrit}
\end{figure}

Finally, we consider meshes with the values of \(\hxm\) and \(\hxp\) being exchanged.
We observe in \Cref{fig:uncrit} a convergence with convenient rates as predicted by \Cref{thm:full}.

\begin{figure}
\centering
\includegraphics[scale=0.2]{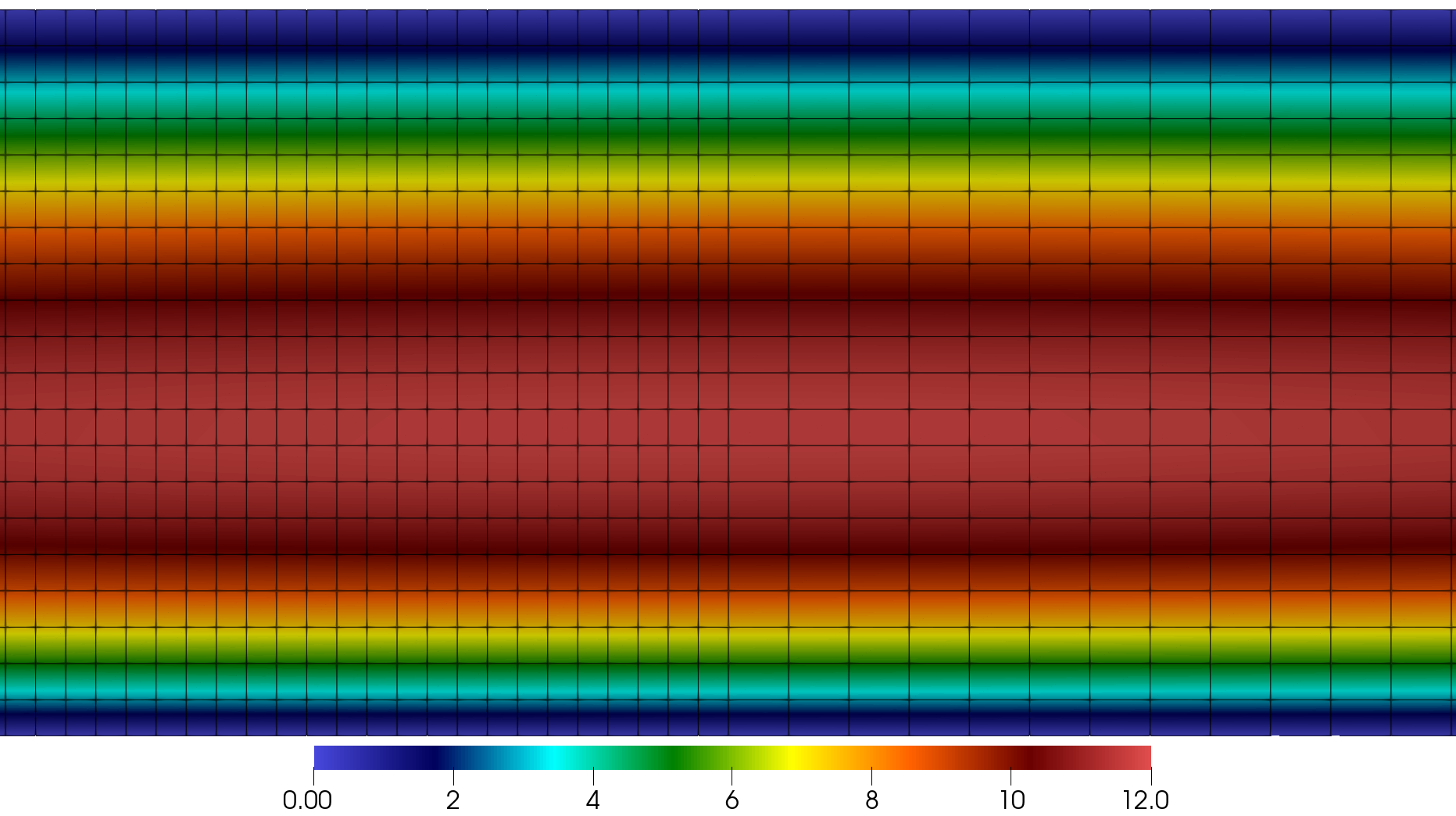}
\caption{
Computed solution for \(\hxp\approx0.26048\) on a flipped mesh in a neighborhood of the interface $\{0\}\times(0,\pi)$.}
\label{fig:uncrit_sol}
\end{figure}
To conclude, for problems with more complicated geometries and discretizations with non-uniform meshes we expect a mixed behaviour, where at each new mesh refinement a stable or unstable setting is dominant in an unpredictable way, giving rise to the commonly observed zigzag error curves.

\bibliographystyle{abbrv}
\bibliography{references}
\end{document}